\documentclass{article}
\usepackage{amsmath,amssymb,euscript,amscd,url,amsthm}
\usepackage{tikz-cd}
\usetikzlibrary{ext.paths.arcto}
\usepackage{tkz-graph}
\usepackage[export]{adjustbox}
\usepackage[title]{appendix}
\usepackage{tabu}
\usepackage{subcaption}
\usepackage{caption}
\usepackage{hyperref}

\title{Combinatorics of double cosets and fundamental domains for the subgroups of the modular group}
\author{Alexey G.~Gorinov, Isaac C.~Kalinkin}
\date{}
\newcommand{\Z}{{\mathbb Z}}
\newcommand{\Q}{{\mathbb Q}}
\newcommand{\C}{{\mathbb C}}
\newcommand{\R}{{\mathbb R}}
\newcommand{\Hh}{{\mathbb H}}
\newcommand{\Id}{\mathrm{Id}}
\newcommand{\ii}{\mathrm{i}}
\newcommand{\rr}{\mathrm{e}^\frac{\pi\mathrm{i}}{3}}
\newcommand{\Xc}{Y_{\mathrm{comb}}}
\newcommand{\PSL}{{\mathrm{PSL}}}
\newcommand{\GL}{{\mathrm{GL}}}
\newcommand{\SL}{{\mathrm{SL}}}
\DeclareMathOperator{\Stab}{{\mathrm{Stab}}}
\DeclareMathOperator{\Aut}{{\mathrm{Aut}}}
\DeclareMathOperator{\Iso}{{\mathrm{Iso}}}

\newtheorem{theorem}{Theorem}[section]
\newtheorem{varthm}{Theorem}

\newtheorem{lemma}[theorem]{Lemma}
\newtheorem{Prop}[theorem]{Proposition}
\newtheorem{cor}[theorem]{Corollary}

\theoremstyle{definition}
\newtheorem{definition}{Definition}[section]
\newtheorem{example}{Example}[section]

\theoremstyle{remark}
\newtheorem{remark}{Remark}[section]
\newtheorem*{question}{Question}

\begin{document}
\maketitle
\begin{abstract}
As noticed by R.~Kulkarni, the conjugacy classes of subgroups of the modular group correspond bijectively to bipartite cuboid graphs. We'll explain how to recover the graph corresponding to a subgroup $G$ of $\mathrm{PSL}_2(\mathbb{Z})$ from the combinatorics of the right action of $\mathrm{PSL}_2(\mathbb{Z})$ on the right cosets $G\backslash\mathrm{PSL}_2(\mathbb{Z})$. This gives a method of constructing nice fundamental domains (which Kulkarni calls "special polygons") for the action of $G$ on the upper half plane. 

For the classical congruence subgroups $\Gamma_0(N)$, $\Gamma_1(N)$, $\Gamma(N)$ etc. the number of operations the method requires is the index times something that grows not faster than a polynomial in $\log N$. This is roughly the square root of the number of operations required by the naive procedure. We give algorithms to locate an element of the upper half-plane on the fundamental domain and to write a given element of $G$ as a product of independent generators. We also (re)prove a few related results about the automorphism groups of modular curves. For example, we give a simple proof that the automorphism group of $X(N)$ is $\mathrm{SL}_2(\mathbb{Z}/N)/\{\pm I\}$.
\end{abstract}
%\section{Introduction}
%\begin{abstract}
%In \cite{kulkarni} R.~Kulkarni gives a method for constructing fundamental domains and finding independent systems of generators for the subgroups of the modular group acting on the upper half-plane. However, the method involves ``trial and error'' and it it not clear how it is supposed to work e.g.\ for $\Gamma(N)$ for general $N$. For classical congruence subgroups $\Gamma_0(N),\Gamma_1(N), \Gamma(N)$ the number of operations our method requires is the index times something that grows more slowly than any positive power of $N$. All proofs are straightforward checks and will only be sketched.
%\end{abstract}

\section{Introduction}
The group $\Gamma=\PSL_2(\Z)$ acts by isometries on the upper half-plane $\Hh=\{z\in\C\mid\mathop{\mathrm{Im}}z>0\}$ equipped with the standard hyperbolic metric. The hyperbolic triangle with vertices $\rr$, $\mathrm{e}^{\frac{2\pi\mathrm{i}}{3}}$ and $\infty$ is a fundamental domain for this action (see e.g.\ \cite[chapter 2]{diashur}). Hence so is the triangle $\triangle$ with vertices $\rr$, $0$ and $\infty$.

\begin{figure}[h!]
\centering
\begin{tikzpicture}
\filldraw[draw=black, thick, fill=gray!20!white] (0,3) -- (0,0) arc(180:120:1) -- (.5,3);
\draw (0,1) arc(90:60:1);
\filldraw[fill=black] (0,1) circle[radius=2pt] ;
\path (-.25,1) node{$\ii$};
\filldraw[fill=white] (.5,.866) circle[radius=2pt] node[anchor=west]{$\rr$};
% $0.8660\approx\frac{\sqrt{3}}{2}$
\draw[->, very thin] (-3,0) -- (3,0);
\path (3,-.3) node{$\mathrm{Re}$};
\path (0,-.3) node{$0$};
\end{tikzpicture}
\caption{Triangle $\triangle$.}\label{fig:0}
\end{figure}

Unless stated otherwise, in the sequel $G$ will be an arbitrary subgroup of $\Gamma$. The union of all $\Gamma$-copies of the hyperbolic line segment in $\Hh$ connecting $\rr$ and $\ii$ is a tree. The image of this tree inside the complex curve $Y(G)=G\backslash\Hh$ is a bipartite cuboid graph, see definition~\ref{def:bipartite_cuboid_graph}. Let us denote this graph $\Xc(G)$. This paper will mainly be concerned with its properties and applications.%It will be one of the main characters in this paper.

Our first observation (see theorem~\ref{main0}) is that one can recover $\Xc(G)$ from the combinatorics of the right action of $\Gamma$ on the set of right cosets $G\backslash\Gamma$. Further, $\Xc(G)$ is a combinatorial model for $Y(G)$, meaning that the geometric realisation of $\Xc(G)$ is homotopy equivalent to $Y(G)$ (theorem~\ref{main0}). These remarks are partially inspired by J.-P.~Serre's book~\cite{serre}.

\smallskip

One of the applications of the graph $\Xc(G)$ is to constructing convex fundamental domains for the action of $G$ on $\Hh$. In \cite{kulkarni} R.~Kulkarni defines special polygons, which are hyperbolic convex polygons of a particular kind plus some combinatorial data, see definition~\ref{special}, and proves that any finite index subgroup of $\Gamma$ admits a special polygon as a fundamental domain, which can be constructed starting from $\Xc(G)$. Moreover, once one has such a fundamental domain, one can easily obtain (\cite[theorem on p.~1055]{kulkarni}) an independent system of generators of $G$, i.e., a subset of $G$ such that $G$ is the free product of the cyclic groups generated by the elements of the subset. In theorem~\ref{main1} we review these results by R.~Kulkarni and remark that the assumption that $G$ should have finite index in $\Gamma$ is unnecessary. We also observe that R.~Kulkarni's theorem on representing $G$ as a free product follows from the results of~\cite{serre}.

R.~Kulkarni shows \cite[theorem 4.2]{kulkarni} that the isomorphism classes of finite bipartite cuboid graphs correspond bijectively to the conjugacy classes of finite index subgroups of $\Gamma$. We prove an analogue of this for subgroups of $\Gamma$ rather than their conjugacy classes. A graph is said to be \emph{pointed} if an edge has been marked on it.

\begin{varthm}[Theorem~\ref{pointed}]
Subgroups of $\Gamma$ correspond bijectively to pointed bipartite cuboid graphs.
\end{varthm}

Again, no finite index assumption is needed.

In remark~\ref{red} we show how one can write a given element of $G$ as a product of the independent generators constructed in section~\ref{gen}; we also explain there how one can ``locate'' a given $z\in\Hh$ on the fundamental domain.

\smallskip

In section~\ref{aut_gr} we review several automorphism groups associated with $G$. The section contains several side results which may be of independent interest. For example, we identify $\Aut(\Xc(G))$ with $N_\Gamma(G)/G$ (proposition~\ref{aut_graph}) and prove that $N_\Gamma(G)/G=N_{\PSL_2(\R)}(G)/G$ provided all cusps of the curve $X(G)$ have the same width (proposition~\ref{equal_width}). We also give a simple proof of the next result:

\begin{varthm}[Proposition~\ref{aut_gamma_n}]
We have $\Aut(X(N))=\SL_2(\Z/N)/\{\pm I\}$.
\end{varthm}

This is Theorem~3.1 in F.~Bars, A.~Kontogeorgis and X.~Xarles~\cite{bkx}, and may have been known earlier. The point of our proof is to illustrate the usefulness of the orbifold Euler characteristic when describing the automorphism groups of modular curves.

\smallskip

Let $N\geq 1$ be an integer. Recall that $\Gamma_0(N)$, $\Gamma_1(N)$, $\Gamma^0(N)$, $\Gamma^1(N)$, and $\Gamma(N)$ are the subgroups of $\Gamma$ such that after lifting to $\SL_2(\Z)$ and reducing modulo $N$ we get the groups of matrices of the form
$$
\begin{pmatrix}a & b\\0& a^{-1}\end{pmatrix},\begin{pmatrix}1 & b\\0& 1\end{pmatrix},\begin{pmatrix}a & 0\\b& a^{-1}\end{pmatrix}, \begin{pmatrix}1 & 0\\b& 1\end{pmatrix},\mbox{ and }\begin{pmatrix}1 & 0\\0& 1\end{pmatrix}
$$
respectively, where $a,b\in\Z/N$ with $a$ invertible. 
%For $G$ equal one of these subgroups the results of section~\ref{grfun} allow one to construct a special polygon much faster than the naive method. Namely, the number of operations is the index (and not the square of the index) times something which grows not faster than a polynomial in $\log N$. 
In section~\ref{classical} we use the results of sections~\ref{grfun} and~\ref{gen} to give algorithms for constructing the graph $\Xc(G)$ and a special polygon for $G$ in each of the series $\Gamma_0(N)$, $\Gamma_1(N)$, $\Gamma^0(N)$, $\Gamma^1(N)$, and $\Gamma(N)$. In~\ref{repg0} and~\ref{repg1g} we show how to construct sets of right coset representatives of $\SL_2(\Z/N)$ modulo the images of $\Gamma_0(N), \Gamma_1(N)$ etc. In section~\ref{constr_the_graph} we give an algorithm for constructing the graph $\Xc(G)$ using the results of~\ref{repg0} and~\ref{repg1g}, and in section~\ref{constr_the_graph} we explain how to construct a fundamental domain for $G$ using $\Xc(G)$. In the next theorem, $O(P\log N)$ stands for ``a sequence $(a_N)_{N\geq 1}$ for which there exists a polynomial $p$ such that $|a_N|<|p(\log N)|$ for all $N\geq 1$''.

\begin{varthm}[Propositions~\ref{constructing_graph} and~\ref{constructing_polygon}]
For $G$ in any of the classical series $\Gamma_0(N), \Gamma^0(N), \Gamma_1(N), \Gamma^1(N)$, or $\Gamma(N)$, the graph $\Xc(G)$ and a fundamental special polygon for the action of $G$ on $\Hh$ can be constructed in $$(\Gamma:G)O(P\log N)$$ operations.
\end{varthm}

This provides an almost quadratic improvement over the naive procedure, which requires $(\Gamma:G)^2$ operations, cf.\ S.-P.~Chan, M.-L.~Lang, C.-H.~Lim and S.-P.~Tan \cite{cllt}.

The algorithms described in sections~\ref{gen} and~\ref{classical} have been implemented on this web page: \url{https://kalinkinisaac.pythonanywhere.com/}. The source files can be found here: \url{https://github.com/kalinkinisaac/modular}.

\begin{remark}
Note that from the computational viewpoint the results of section~\ref{classical} can only be improved by a factor of a polynomial in $\log N$ as the output itself is of the order $(\Gamma:G)$. For example, the number of edges in a special polygon is twice the number of independent generators of $G$, and the latter is $\in \left(\frac{1}{6}(\Gamma:G),\frac{1}{3}(\Gamma:G)+3\right)$, as one can see e.g.\ by using the formula of \cite[Theorem 3.1.1]{diashur} and noting that in this formula $d=(\Gamma:G)$, the number of independent generators is $\varepsilon_2+\varepsilon_3+2g+\varepsilon_\infty-1$, and $\varepsilon_\infty>0, \varepsilon_2,\varepsilon_3, g\geq 0$.
\end{remark}

\begin{remark}
In~\cite{nie_parent} Z.~Nie and C.~X.~Parent give another procedure for constructing fundamental domains for the classical congruence subgroups. It has the advantage that the coset representatives are given by an explicit formula. The advantage of our approach is that the fundamental domain one obtains is automatically convex. For comparison see e.g.\ figure~\ref{fig:spec_poly_gamma_6} in the present paper and Figure~1 in~\cite{nie_parent}, which both show fundamental domains for $\Gamma_0(6)$. Z.~Nie and C.~X.~Parent do not give a complexity estimate for their algorithm, and we were not able to compare their method with ours in terms of complexity.
\end{remark}

%The initial motivation for this work was to find a clear and fast procedure for constructing a special polygon for $\Gamma(N)$; these polygons are needed to give an explicit description of the mixed Hodge structure on the cohomology of the symmetric powers of the Hodge local system on $X(N)$, see~\cite{elliptic}.

\begin{question}
If we take a bipartite cuboid graph, forget the partition of the vertices into the two types, and replace any two edges meeting at a vertex of valency 2 by a single edge, we will obtain a unitrivalent graph with a cyclic orientation at each trivalent vertex. Such graphs occur in the theory of Vassiliev (aka finite type) invariants of knots in 3-space (see e.g.~\cite{bn}). A natural question is whether this theory is related in some way to subgroups of $\Gamma$. In particular, are there any natural functions on the set of the conjugacy classes of finite index subgroups of $\Gamma$ satisfying the anti-symmetry and IHX relations?
\end{question}

{\bf Notation and conventions.} We will need to distinguish between left and right actions and quotients. If $X$ is a set and a group $H_1$, respectively a group $H_2$ acts on $X$ on the left, respectively on the right, the quotients will be denoted $H_1\backslash X$, respectively $X/H_2$. If in addition the actions commute, we have the double quotient $H_1\backslash X/H_2=(H_1\backslash X)/H_2=H_1\backslash (X/H_2)$.

Suppose $x,y$ are two distinct points in $\bar\Hh=\Hh\sqcup\R\sqcup\{\infty\}$. There is a unique geodesic $\gamma\subset\Hh$ such that the closure $\bar\gamma\subset\bar\Hh$ contains $x$ and $y$. We use $[x,y], (x,y), [x,y), (x,y]$ to denote the corresponding closed, open and half-open intervals in $\bar\gamma$.

A {\it cyclic ordering} of a set $X$ with $n$ elements is a free action of the cyclic group $\Z/n$ on $X$. If such an action is given, we say that $1\cdot x$, respectively $(-1)\cdot x$ {\it comes after} $x\in X$, respectively {\it comes before} $x$ in the cyclic ordering. In this paper we will only consider cyclic orderings of three element sets.

\smallskip

{\bf Acknowledgements.} The first named author is grateful to V.~Gritsenko, A.~Kalmynin, S.~Leli\`evre, and H.~Verrill for helpful and interesting discussions. The first version of the paper was written while the first named author was visiting the Max Planck Institute for Mathemathics in Bonn. This author is grateful to MPIM Bonn for excellent working conditions and stimulating atmosphere.

The authors are grateful to the anonymous reviewer for carefully reading the paper and making useful suggestions.

\section{Graphs and fundamental domains}\label{grfun}
%In the sequel we assume $(\Gamma:G)<\infty$ to ensure that all graphs and polygons we consider are finite (respectively, are made of a finite number of $\Gamma$-copies of some standard polygon).

We will write $\left[\begin{smallmatrix}a&b\\c&d\end{smallmatrix}\right]$ to denote the image of the matrix $\left(\begin{smallmatrix}a&b\\c&d\end{smallmatrix}\right)\in\mathrm{SL}_2(\Z)$ in $\Gamma$. 

Consider the set $G\backslash\Gamma$, on which $\Gamma$ acts on the right. Set $G_0$ and $G_1$ to be the subgroups of $\Gamma$ respectively of order 2 and 3 generated respectively by $\left[\begin{smallmatrix}0&-1\\1&0\end{smallmatrix}\right]$ and $\left[\begin{smallmatrix}0&1\\-1&1\end{smallmatrix}\right]$.

Let $\Xc(G)$ be the graph constructed as follows. The set $\mathop{\mathrm{Vert}}(\Xc(G))$ of vertices of $\Xc(G)$ is the disjoint union of the sets of all $G_0$- and $G_1$-orbits in $G\backslash\Gamma$, i.e.\ $$\mathop{\mathrm{Vert}}(\Xc(G))=(G\backslash\Gamma/G_0)\sqcup (G\backslash\Gamma/G_1).$$ The set $\mathop{\mathrm{Edge}}(\Xc(G))$ of edges of $\Xc(G)$ is $G\backslash\Gamma$. Viewing the elements of $\mathop{\mathrm{Vert}}(\Xc(G))$ as subsets of $G\backslash \Gamma$, we declare that an element $x\in\mathop{\mathrm{Edge}}(\Xc(G))=G\backslash\Gamma$ joins the vertices $A$ and $B$, $A,B\subset G\backslash\Gamma$ iff $x\in A\cap B$.
%Note that if two vertices of $\Xc(G)$ are connected with an edge, then one is a $G_0$-orbit and the other one is a $G_1$-orbit.) The graph $\Xc(G)$ has a ``distinguished'' edge $\mathbf{r}_G$ which is $G\in G\backslash\Gamma$ (i.e.\ $G$ as a right $G$-coset). The edge $\mathbf{r}_G$ joins two ``distinguished'' vertices, which are the $G_0$- and $G_1$-orbits of $G\in G\backslash\Gamma$; denote these vertices respectively by $V^0_G$ and $V^1_G$.
Note that the vertices of $\Xc(G)$ can have valency 1, 2 or 3.

Let us recall the following definition (see~\cite[4.1]{kulkarni}).

\begin{definition}\label{def:bipartite_cuboid_graph}
A graph $X$ such that
\begin{itemize}\label{graph}
\item the vertices of $X$ are subdivided into two types;
\item the vertices of one type can have valency 1 or 2;
\item the vertices of the other type can have valency 1 or 3;
\item any edge of $X$ joins a vertex of one type with a vertex of the other type;
\item there is a cyclic order on the set of the edges meeting at a vertex of valency~3
\end{itemize}
is called a {\it bipartite cuboid graph}. 
\end{definition}

Observe that the graph $\Xc(G)$ has a bipartite cuboid structure. We will say that the vertices of $\Xc(G)$ which are elements of $G\backslash\Gamma/G_0$ are \emph{of type~(0)}, and those which are elements of $G\backslash\Gamma/G_1$ are \emph{of type~(1)}. The vertices of type~(0) can be of valency 2 or 1, and the vertices of type~(1) can be of valency 3 or 1. Any edge of $\Xc(G)$ joins a vertex of type~(0) and a vertex of type~(1). The edges of $\Xc(G)$ which meet at a vertex of valency 3 (and type~(1)) are cyclically ordered: they correspond to the elements of a 3-element $G_1$-orbit $A$; take any element $x\in A$ and the cyclic order will be $x,x\cdot \left[\begin{smallmatrix}0&1\\-1&1\end{smallmatrix}\right],x\cdot\left[\begin{smallmatrix}0&1\\-1&1\end{smallmatrix}\right]^2.$

\smallskip

Before going further, let us consider some examples. In figure~\ref{fig:1} below the vertices of type~(0) are shown in black, and the vertices of type~(1) in white. The cyclic orientation is induced from the standard orientation of the plane: the group $\Z/3$ then acts counter-clockwise on the edges that meet at any trivalent vertex of type~(1).

\begin{example}\label{example_graphs}
(a) $G=\Gamma$. In this case $G\backslash\Gamma$ contains just one element, so the graph $\Xc(G)$ is as shown in figure~\ref{fig:1a}.

\begin{figure}[h!]
\centering
\begin{minipage}{0.3\textwidth}
  \centering
\begin{tikzpicture}
\SetGraphUnit{2.5}
\SetVertexNoLabel
\tikzset{VertexStyle/.style = {
shape = circle,
fill = black,
inner sep = 0pt,
outer sep = 0pt,
minimum size = 5pt,
draw}}
\Vertex{A}
\tikzset{VertexStyle/.style = {
shape = circle,
fill = white,
inner sep = 0pt,
outer sep = 0pt,
minimum size = 5pt,
draw}}
\EA[unit=2](A){B}
\Edge(A)(B)
\end{tikzpicture}
\subcaption{$\Xc(\Gamma)$.}\label{fig:1a}
\end{minipage}
\begin{minipage}{0.3\textwidth}
\centering
\begin{tikzpicture}
\SetGraphUnit{2.5}
\SetVertexNoLabel
\tikzset{VertexStyle/.style = {
shape = circle,
fill = black,
inner sep = 0pt,
outer sep = 0pt,
minimum size = 5pt,
draw}}
\Vertex{A}
\SO[unit=2](A){B}
\SO[unit=4](A){C}
\tikzset{VertexStyle/.style = {
shape = circle,
fill = white,
inner sep = 0pt,
outer sep = 0pt,
minimum size = 5pt,
draw}}
\WE[unit=2](B){D}
\EA[unit=2](B){E}
\Edge(B)(E)
\Edge(B)(D)
\Edge[style={bend left=45}](D)(A)
\Edge[style={bend right=45}](D)(C)
\Edge[style={bend right=45}](E)(A)
\Edge[style={bend left=45}](E)(C)
\end{tikzpicture}
\subcaption{$\Xc(\Gamma(2))$.}\label{fig:1b}
\end{minipage}
\begin{minipage}{0.3\textwidth}
  \centering
\begin{tikzpicture}
\SetGraphUnit{2.5}
\SetVertexNoLabel
\tikzset{VertexStyle/.style = {
shape = circle,
fill = black,
inner sep = 0pt,
outer sep = 0pt,
minimum size = 5pt,
draw}}
\Vertex{A}
\EA[unit=4](A){C}
\tikzset{VertexStyle/.style = {
shape = circle,
fill = white,
inner sep = 0pt,
outer sep = 0pt,
minimum size = 5pt,
draw}}
\EA[unit=2](A){B}
\Edge(B)(C)
\Edge[style={bend left=90},label=$\mathbf{r}_{\Gamma_0(2)}$, labelstyle={below}](B)(A)
\Edge[style={bend right=90},label=$\mathbf{r}_{\Gamma^0(2)}$, labelstyle={above}](B)(A)
\end{tikzpicture}
\subcaption{$\Xc(\Gamma_0(2))=\Xc(\Gamma^0(2))$ and the distinguished edges, see definition~\ref{dist_edge_def}.}\label{fig:1c}
\end{minipage}
\begin{minipage}{0.3\textwidth}
\centering
\begin{tikzpicture}
\SetGraphUnit{2.5}
\SetVertexNoLabel
\tikzset{VertexStyle/.style = {
shape = circle,
fill = black,
inner sep = 0pt,
outer sep = 0pt,
minimum size = 5pt,
draw}}
\Vertex{A}
\tikzset{VertexStyle/.style = {
shape = circle,
fill = white,
inner sep = 0pt,
outer sep = 0pt,
minimum size = 5pt,
draw}}
\EA[unit=2](A){B}
\WE[unit=2](A){C}
\Edge(A)(B)
\Edge(A)(C)
\end{tikzpicture}
\subcaption{$\Xc(\ker(\Gamma\to\Z/2))$.}\label{fig:1d}
\end{minipage}
\begin{minipage}{0.3\textwidth}
\begin{tikzpicture}
\SetGraphUnit{2.5}
\SetVertexNoLabel
\tikzset{VertexStyle/.style = {
shape = circle,
fill = white,
inner sep = 0pt,
outer sep = 0pt,
minimum size = 5pt,
draw}}
\Vertex{A}
\tikzset{VertexStyle/.style = {
shape = circle,
fill = black,
inner sep = 0pt,
outer sep = 0pt,
minimum size = 5pt,
draw}}
\Vertex[a=60 , d=2 cm]{B}
\Vertex[a=180 , d=2 cm]{C}
\Vertex[a=-60 , d=2 cm]{D}
\Edge(A)(B)
\Edge(A)(C)
\Edge(A)(D)
\end{tikzpicture}
\subcaption{$\Xc(\ker(\Gamma\to\Z/3))$.}\label{fig:1e}
\end{minipage}
\caption{Examples of $\Xc(G)$. The vertices in $G\backslash\Gamma/G_0$, respectively in $G\backslash\Gamma/G_1$ are shown in black, respectively in white.} \label{fig:1}
\end{figure}

(b) $G=\Gamma(2)$. The subgroup $G$ is normal in $\Gamma$, and the quotient $G\backslash\Gamma$ is isomorphic to $\PSL_3(\Z/2)\cong S_3$. The right action of $\Gamma$ on $G\backslash\Gamma$ is isomorphic to the action of $\Gamma$ on $S_3$ by right shifts via the homomorphism $f:\Gamma\to S_3$ that takes $\left[\begin{smallmatrix}0&-1\\1&0\end{smallmatrix}\right]$ to $(12)$ and $\left[\begin{smallmatrix}0&1\\-1&1\end{smallmatrix}\right]$ to $(123)$. The structure of $G_0$- and $G_1$-orbits for this action is shown in figure~\ref{fig:2a}, so the graph $\Xc(G)$ is as shown in figure~\ref{fig:1b}.

\begin{figure}[h!]
\centering
\begin{minipage}{0.3\textwidth}
\centering
\begin{tikzpicture}
\draw (0,-.5) ellipse(50pt and 10pt);
\draw (0,.5) ellipse(50pt and 10pt);
\fill[black] (0,-.5) circle(3pt);
\fill[black] (.7,-.5) circle(3pt);
\fill[black] (-.7,-.5) circle(3pt);
\fill[black] (0,.5) circle(3pt);
\fill[black] (.7,.5) circle(3pt);
\fill[black] (-.7,.5) circle(3pt);
\draw (-.3,-.7) rectangle(.3,.7);
\draw (-1,-.7) rectangle(-.4,.7);
\draw (.4,-.7) rectangle(1,.7);
\end{tikzpicture}
\subcaption{$G=\Gamma(2)$.}\label{fig:2a}
\end{minipage}
\begin{minipage}{0.3\textwidth}
\centering
\begin{tikzpicture}
\draw (0,0) ellipse(50pt and 20 pt);
\fill[black] (0,0) circle(3pt);
\fill[black] (.7,0) circle(3pt);
\fill[black] (-.7,0) circle(3pt);
\draw (-.3,-.3) rectangle(1,.3);
\draw (-1,-.3) rectangle(-.4,.3);
\end{tikzpicture}
\subcaption{$G=\Gamma_0(2)$.}\label{fig:2b}
\end{minipage}
\begin{minipage}{0.3\textwidth}
\centering
\begin{tikzpicture}
\draw (-1.75,-.6) rectangle(1.75,.6);
\fill[black] (-.7,0) circle(3pt);
\fill[black] (.7,0) circle(3pt);
\draw (-.7,0) circle(10pt);
\draw (.7,0) circle(10pt);
\end{tikzpicture}
\subcaption{$G=\ker(\Gamma\to\Z/2)$.}\label{fig:2c}
\end{minipage}
\caption{Orbits of $G_0$ and $G_1$ in $G\backslash\Gamma$. The rectangular boxes represent the $G_0$-orbits, and the elliptic ones the $G_1$-orbits.}\label{fig:2}
\end{figure}

(c) Since $\Gamma(2)\triangleleft\Gamma$, the group $\Gamma$ acts on $\Gamma(2)\backslash\Gamma$ not only on the right but also on the left. This is the action on $S_3=\Gamma(2)\backslash\Gamma$ by left shifts via the homomorphism $f$ from the previous example. If $G\subset\Gamma$ is a subgroup that contains $\Gamma(2)$, then the right $\Gamma$-set $G\backslash\Gamma$ is isomorphic to the quotient of $\Gamma(2)\backslash \Gamma$ by the left action of $G\subset\Gamma$, or equivalently to $f(G)\backslash S_3$ with the right action of $\Gamma$ induced by $f$. 

Now take $G=\Gamma_0(2)=\Gamma_1(2)$. Then the group $f(G)$ will be generated by $$f\left(\left[\begin{smallmatrix}1&-1\\0&1\end{smallmatrix}\right]\right)=f\left(\left[\begin{smallmatrix}0&-1\\1&0\end{smallmatrix}\right]\right)f\left(\left[\begin{smallmatrix}0&1\\-1&1\end{smallmatrix}\right]\right)=(12)(123)=(23).$$ (Here and in the sequel we multiply permutations by first applying the one on the right.) The $G_0$- and $G_1$-orbits for the (right) action of $\Gamma$ on $G\backslash\Gamma=f(G)\backslash S_3$ are shown in figure~\ref{fig:2b}, so we get the graph $\Xc(G)$, see figure~\ref{fig:1c}.

If we take $G=\Gamma^0(2)=\Gamma^1(2)$, then one can check as above that $f(G)$ is generated by $(13)=(123)(12)$, and that the graph $\Xc(G)$ is the same as for $\Gamma_0(2)$. (Both graphs have distinguished edges, see definition~\ref{dist_edge_def}; these are different for $\Gamma_0(2)$ and $\Gamma^0(2)$ and are also shown in figure~\ref{fig:1c}.)

(d) Arguing as in the previous example but taking $G=f^{-1}(A_3)$, we get the graph $\Xc(G)$ shown in figure~\ref{fig:1d} using the orbit structure for the actions of $G_0$ and $G_1$ shown in figure~\ref{fig:2c}. Note that in this case $G$ is the kernel of the unique non-trivial group homomorphism $\Gamma\to\Z/2$.

(e) Finally, let $G$ be the kernel of a non-trivial group homomorphism $g:\Gamma\to\Z/3$. (There are two such homomorphisms, but their kernels are the same.) Observe that $g$ takes $\left[\begin{smallmatrix}0&1\\-1&1\end{smallmatrix}\right]$ to a generator of $\Z/3$ and $\left[\begin{smallmatrix}0&-1\\1&0\end{smallmatrix}\right]$ to $0$. Arguing as above we see that the graph $\Xc(G)$ is as shown in figure~\ref{fig:1e}.

\end{example}

In the sequel we will need another description of the graph $\Xc(G)$. Set $\mathbf{T}=\Xc(\{\Id\})$. The graph $\mathbf{T}$ is an infinite tree with $\mathop{\mathrm{Edge}}(\mathbf{T})=\Gamma$, $$\mathop{\mathrm{Vert}}(\mathbf{T})=(\Gamma/G_0)\sqcup(\Gamma/G_1),$$ and $g\in\Gamma$ joins $g_0G_0$ and $g_1G_1$ iff $g\in g_0G_0\cap g_1G_1$. 

The group $\Gamma$ acts on $\mathbf{T}$ on the left, and we have an isomorphism $G\backslash\mathbf{T}=\Xc(G)$.

\begin{figure}[h!]
\centering
\begin{tikzpicture}
\SetGraphUnit{2.5}
\SetVertexNoLabel
\tikzset{VertexStyle/.style = {
shape = circle,
fill = white,
inner sep = 0pt,
outer sep = 0pt,
minimum size = 5pt,
draw}}
\Vertex{A0}
\Vertex[a=60 , d=2 cm]{A1}
\Vertex[a=180 , d=2 cm]{A2}
\Vertex[a=-60 , d=2 cm]{A3}
\tikzset{VertexStyle/.style = {
shape = circle,
fill = black,
inner sep = 0pt,
outer sep = 0pt,
minimum size = 5pt,
draw}}
\Vertex[a=60 , d=1 cm]{B1}
\Vertex[a=180 , d=1 cm]{B2}
\Vertex[a=-60 , d=1 cm]{B3}
\Vertex[x=2,y=1.7321]{B11}
\Vertex[x=.5,y=2.5981]{B12}
\Vertex[x=-2.5,y=.866]{B21}
\Vertex[x=-2.5,y=-.866]{B22}
\Vertex[x=2,y=-1.7321]{B31}
\Vertex[x=.5,y=-2.5981]{B32}
\tikzset{VertexStyle/.style = {}}
\Vertex[x=3,y=1.7321]{C11}
\Vertex[x=0,y=3.4641]{C12}
\Vertex[x=-3,y=1.7321]{C21}
\Vertex[x=-3,y=-1.7321]{C22}
\Vertex[x=3,y=-1.7321]{C31}
\Vertex[x=0,y=-3.4641]{C32}
\Edge(A0)(B1)
\Edge(A0)(B2)
\Edge(A0)(B3)
\Edge(A1)(B1)
\Edge(A2)(B2)
\Edge(A3)(B3)
\Edge(A1)(B11)
\Edge(A1)(B12)
\Edge(A2)(B21)
\Edge(A2)(B22)
\Edge(A3)(B31)
\Edge(A3)(B32)
\tikzset{EdgeStyle/.style = dotted}
\Edge(B11)(C11)
\Edge(B12)(C12)
\Edge(B21)(C21)
\Edge(B22)(C22)
\Edge(B31)(C31)
\Edge(B32)(C32)
\end{tikzpicture}
\caption{Tree $\mathbf{T}=\Xc(\{\Id\})$.}\label{fig:3}
\end{figure}

Let us now relate the graph $\Xc(G)$ to the complex curve $Y(G)=G\backslash\Hh$. Set $\mathcal{T}$ to be the union of all $\Gamma$-copies of the hyperbolic line segment $\mathbf{r}=[\ii,\rr]$. It is not difficult to show that $\mathcal{T}$ is a tree\footnote{Here and in the sequel we do not distinguish between a graph and its geometric realisation; this should never lead to a confusion.}. Since both quotient graphs $\Gamma\backslash\mathcal{T}$ and $\Gamma\backslash\mathbf{T}$ are line segments, it follows from \cite[I.4, Theorem 7]{serre} (or can be easily checked by hand) that there is a $\Gamma$-equivariant isomorphism $\mathbf{T}\to\mathcal{T}$. Let $f:\mathbf{T}\to\mathcal{T}$ be the unique $\Gamma$-equivariant isomorphism such that the vertex of $\mathbf{T}=\Xc(\{\Id\})$ that corresponds to $G_0$ as a left $G_0$-coset goes to $\ii$, and the vertex that corresponds to $G_1$ as a left $G_1$-coset goes to $\rr$.

Recall that we have set $\triangle$ to be the hyperbolic triangle with vertices $0,\rr,\infty$. %minus the vertices $0$ and $\infty$.
%It is well known that $\triangle$ is a fundamental domain for the action of $\Gamma$ on $\mathbb{H}$. 
%Let $p:\triangle\to\mathbf{r}$ be a strong deformation retraction that takes the imaginary half-axis to $\ii$ and both the geodesic intervals $[\rr,0)$ and $[\rr,\infty)$ to $\rr$. Extend $p$ to 
Let $p:\mathbb{H}\to\mathcal{T}$ be a $\Gamma$-equivariant strong deformation retraction. %The mapping $p$ can be joined with $\mathop{\Id_\mathbb{H}}$ by a $\Gamma$-equivariant homotopy leaving $\mathcal{T}$ pointwise fixed. 
From the above we obtain the following

\begin{theorem}\label{main0}
\begin{enumerate}
\item The graph $\Xc(G)$ is a combinatorial model for the Riemann surface $Y(G)=G\backslash\mathbb{H}$, i.e.\ $\Xc(G)$ is homotopy equivalent to $Y(G)$.
\item Cut $\Xc(G)$ in some of the vertices of valency 2 to obtain a tree $Y'$ and take a lift $\jmath$ of this tree into $\mathbf{T}$ so that the diagram
\begin{equation}\label{diag}
\begin{tikzcd}
Y'\ar[r,"\jmath"] \ar[dr] &\mathbf{T}\ar[d]\\
&\Xc(G)=G\backslash\mathbf{T}
\end{tikzcd}
\end{equation}
commutes. The set $P=p^{-1}\big(f(\jmath(Y'))\big)$ is a convex hyperbolic polygon which is a union of $\Gamma$-copies of $\triangle$, and a fundamental domain for the action of $G$ on $\mathbb{H}$.
\end{enumerate}
\end{theorem}

{\bf Proof.} Since the maps $p$ and $f$ are $\Gamma$-equivariant, the quotient $G\backslash\Hh$ is homotopy equivalent to $G\backslash\mathcal{T}$, and the latter is homeomorphic to $G\backslash\mathbf{T}=\Xc(G)$. This proves the first part of the theorem. To prove the second part observe that $\jmath(Y')$ is a fundamental domain for the action of $G$ on $\mathbf{T}$, so $f(\jmath(Y'))$ is a fundamental domain for the action of $G$ on $\mathcal{T}$.$\clubsuit$

\smallskip

The graph $\Xc(G)$ has a ``distinguished'' edge:

\begin{definition}\label{dist_edge_def} The \emph{distinguished edge} of edge of $\Xc(G)$, denoted $\mathbf{r}_G$, is $G\in G\backslash\Gamma$ (i.e.\ $G$ as a right $G$-coset). The edge $\mathbf{r}_G$ joins two \emph{distinguished vertices}, which are the $G_0$- and $G_1$-orbits of $G\in G\backslash\Gamma$; we denote these vertices respectively by $V_0^G$ and $V_1^G$. 
\end{definition}

Note that the distinguished edge of $\mathbf{T}=\Xc(\{\Id\})$ is the identity $\Id\in\Gamma=\{\Id\}\backslash\Gamma$. The isomorphism $f:\mathbf{T}\to\mathcal{T}$ was defined so that the distinguished edge goes to $\mathbf{r}=[\ii,\rr]$.

\smallskip

\begin{example} 
For all graphs shown in figure~\ref{fig:1}, except $\Xc(\Gamma_0(2))$, the group of automorphisms (that preserve the types of the vertices and the cyclic orderings at each trivalent vertex of type~(1)) is transitive on the edges. As we will see later (theorem~\ref{pointed}) this reflects the fact that the corresponding subgroups are normal in $\Gamma$. So one can view any edge as distinguished in these cases. 

Now take $G=\Gamma_0(2)$. As we saw earlier, $G\backslash\Gamma$ is isomorphic to $H\backslash S_3$ where $H\subset S_3$ is generated by $(23)$, and $\left[\begin{smallmatrix}0&-1\\1&0\end{smallmatrix}\right]$, respectively $\left[\begin{smallmatrix}0&1\\-1&1\end{smallmatrix}\right]$ acts by right multiplication by $(12)$, respectively by $(123)$. Observe that
$$\{(13),(123)\}\cdot (12)=\{(13),(123)\}, \mbox{ and } \{(13),(123)\}\cdot (123)^{-1}=\{e, (23)\}.$$
This means that (a) the edge $\{(13),(123)\}\in H\backslash S_3=G\backslash\Gamma$ of $\Xc(G)$ connects the unique vertex of type~(1) with the unique univalent vertex of type~(0), and (b) the distinguished edge $\{e,(23)\}$ comes before $\{(13),(123)\}$ is the cyclic ordering. So we see that the distinguished edge of $\Xc(G)$ is as shown in figure~\ref{fig:1c}.

A similar calculation for $\Gamma^0(2)$ results in the same graph but a different distinguished edge, also shown in figure~\ref{fig:1c}.
\end{example}

Now let $Y'$ and $\jmath$ be as in theorem~\ref{main0}, and let $\imath:Y'\to \mathcal{T}$ be $f\circ\jmath$.
\iffalse
Let $\mathbf{r}'$ be the edge of $Y'$ that corresponds to $\mathbf{r}_G$. %, and let $V'_0$ and $V'_1$ be the end points of $\mathbf{r}'$; these map to $V_0^G$, respectively $V_1^G$.
Observe that the tree $\mathcal{T}$ is bipartite cuboid as it is isomorphic to $\mathbf{T}=\Xc(\{\Id\})$; the vertices of type~(0), respectively of type~(1) are the $\Gamma$-images of $\ii$, respectively of $\rr$, and the cyclic orientation at each vertex of type~(1) is induced by the counter-clockwise rotation through $\frac{2\pi}{3}$ about the vertex.

Let $\imath:Y'\to \mathcal{T}$ be $f\circ\jmath$. The map $\imath$ takes $\mathbf{r}'$ to $[\ii,\rr]$ and preserves the types of the vertices and the cyclic orientation at each trivalent vertex. 
\fi
Each edge $r$ of $\mathcal{T}$ can be written as $r=g_r\cdot\mathbf{r}$ for a unique $g_r\in\Gamma$.
\begin{Prop}\label{fund_domain_from_graph}
The polygon $P$ of theorem \ref{main0} is the union of $g_r\cdot\triangle$ for all edges $r$ of $\mathcal{T}$ which are in the image of $\imath$.
\end{Prop}

{\bf Proof.} This follows from the fact that each edge $s$ of $\jmath(Y')$ is the image of the distinguished edge of $\mathbf{T}$ under a unique $g_s\in\Gamma$, and $f$ takes the distinguished edge of $\mathbf{T}$ to $\mathbf{r}$.$\clubsuit$

\smallskip

%Hindsight: I wish I had put this in the journal version.
\begin{remark}\label{rmk:what_distinguished_edges_are_needed_for}
The tree $Y'$ as in theorem~\ref{main0} inherits a distinguished edge from $\Xc(G)$. There is a unique embedding $\jmath:Y'\to\mathbf{T}$ which takes the distinguished edge to the distinguished edge. Moreover, the diagram~(\ref{diag}) then commutes. This allows us to construct the fundamental polygon $P$ once we know the graph $\Xc(G)$ and its distinguished edge.
\end{remark}
%End of hindsight

\section{Special polygons and independent sets of generators}\label{gen}

Let us recall the following definition from~\cite[2.4]{kulkarni}:

\begin{definition}\label{special}
A {\it special polygon} is a convex hyperbolic polygon $P\subset\mathbb{H}$ plus an involution $\sigma$ on the set of the edges of $P$ which
satisfy the following.

\begin{enumerate}
\item Every edge of $P$ is a $\Gamma$-copy of the geodesic interval that joins $\infty$ with either $0,\ii$, or $\rr$.
\item Every edge which is a $\Gamma$-copy of the geodesic $(0,\infty)$ is paired under~$\sigma$ with another such edge.
\item Every edge $s$ which is a $\Gamma$-copy of the geodesic interval $[\ii, \infty)$ is adjacent to another such interval $s'$, and $\sigma(s)=s',\sigma(s')=s$. The angle between $s$ and $s'$ is $\pi$.
\item Every edge $s$ which is a $\Gamma$-copy of the geodesic interval $[\rr,\infty)$ is adjacent to another such line segment $s'$, and $\sigma(s)=s',\sigma(s')=s$. The internal angle between $s$ and $s'$ is $\frac{2\pi}{3}$.
\end{enumerate}
\end{definition}

Let now $P, Y', \jmath:Y'\to\mathbf{T}, \imath: Y'\to\mathcal{T}$ be as in the previous section. Let us show that the polygon $P$ is special. To do this we need to define an involution $\sigma$ on the set of edges of $P$, which correspond to the univalent vertices of $Y'$. More precisely:

\begin{itemize}\label{edges}
\item A univalent vertex of type~(0) can be obtained either from cutting $\Xc(G)$ at a bivalent vertex, or from a univalent vertex of $\Xc(G)$. In the first case the vertex gives us a single edge, and in the second case two edges such that the angle between them is $\pi$. 
\item For each univalent vertex of type~(1) there are two edges of $P$, and the internal angle between them is~$\frac{2\pi}{3}$.
\end{itemize}

We will now define the involution $\sigma$. Moreover, for each edge $s$ of $P$ we will define an element $g_s\in G$ that takes $s$ to $\sigma(s)$; we will shortly need these elements of $G$ in order to state theorem~\ref{main1}.

Suppose an edge $s$ of $P$ corresponds to a univalent vertex $V$ of $Y'$ obtained from cutting a bivalent vertex $W$ of $\Xc(G)$ of type~(0), and let $\bar{s}$ be the edge of $P$ corresponding to the other vertex $\bar{V}$ of $Y'$ obtained from cutting $W$. Then set $\sigma(s)=\bar{s},\sigma(\bar{s})=s$. Both $s$ and $\bar{s}$ are $\Gamma$-copies of the imaginary half-axis. Equip $s$ and $\bar{s}$ with the orientation induced from $P$; there is a unique $g_s\in\Gamma$ that takes $s$ to $\bar{s}$ in an orientation-reversing way. Moreover, $g_s$ takes $V$ to $\bar{V}$ (recall that we have identified $Y'$ with $\jmath(Y')$) and $Y'\cap g_s(Y')=\{\bar{V}\}$; this implies that $g_s\in G$ (since $Y'$ is a fundamental domain for the action of $G$ on $\mathbf{T}$).

Suppose edges $s', s''$ correspond to a univalent vertex $V$ of $\Xc(G)$ of type~(0). There is a $g\in \Gamma$ such that $s'=g([\ii,0))$, and $s''=g([\ii,\infty))$. Define $\sigma(s')=s'',\sigma(s'')=s$. There is a unique element $g'\in\Gamma$ that swaps $s'$ with $s''$. Note that $g'^2=\Id$; set $g_{s'}=g_{s''}=g'$. Moreover, $g'\cdot V=V$ and $Y'\cap (g'\cdot Y')=\{V\}$, which implies that $g'\in G$.

Finally, suppose we have two edges $s'$ and $s''$ of $P$ which correspond to a univalent vertex $V$ of $Y'$ obtained from a univalent vertex of $\Xc(G)$ of type~(1). Then there is a unique $g\in\Gamma$ such that $s'$ and $s''$ are the images under $g$ of the geodesic intervals $[\rr,0),[\rr,\infty)\subset\mathbb{H}$. Set $\sigma(s')=s'',\sigma(s'')=s$. There is a unique $g'\in\Gamma$ that takes $s'$ to $s''$; set $g_{s'}=g',g_{s''}=g'^{-1}=g'^2$. As above we have $g_{s'},g_{s''}\in G$.

\begin{remark} 
The elliptic points of $Y(G)$ can be located on the fundamental domain $P$ as follows: the elliptic points of order 2 are the points of the form $\imath(V)$ where $V$ is a univalent vertex of $Y'$ of type~(0) that comes from a vertex of $\Xc(G)$ of type~(0); the elliptic points of order 3 are those that have the form $\imath(V)$, $V$ a univalent vertex of $Y'$ of type~(1).
\end{remark}

\begin{example}
Let us describe $P$ and $\sigma$ for the subgroups from example~\ref{example_graphs}. Note that in each of these examples there is a unique way (up to automorphism) to cut $\Xc(G)$ at vertices of type~(0) to get the tree $Y'$. We also need to specify the embedding $\imath:Y'\to\mathcal{T}$. We do this by requiring that the distinguished edge should go to $\mathbf{r}=[\ii,\rr]$.

In figure~\ref{fig:example_graphs_contd} the polygon $P$ is the shaded area, and the tree $\imath(Y')$ is drawn inside it. The vertices of type~(0), respectively of type~(1) are marked $\bullet$, respectively $\circ$. Pairs of edges which contain $\bullet$'s with the same label are interchanged by $\sigma$. If an edge contains an unlabelled $\bullet$, we regard it as the union of two edges, which are paired under $\sigma$. The same applies to the pair of edges meeting at a $\circ$. The elliptic points of order $2$, respectively $3$ are the unlabelled $\bullet$'s, respectively the $\circ$'s on the boundary of $P$. The cusps are the vertices on the boundary of $\Hh$, i.e.\ in $\R\sqcup\{\infty\}$.

\begin{figure}[h!]
\centering
\begin{minipage}{0.3\textwidth}
\centering
\begin{tikzpicture}
\filldraw[draw=black, thick, fill=gray!20!white] (0,3) -- (0,0) arc(180:120:1) -- (.5,3);
\draw (0,1) arc(90:60:1);
\filldraw[fill=black] (0,1) circle[radius=2pt] ;
\filldraw[fill=white] (.5,.866) circle[radius=2pt];
% $0.866\approx\frac{\sqrt{3}}{2}$
\path (0,-.3) node{$0$};
\end{tikzpicture}
\subcaption{$G=\Gamma$.}\label{fig:p_for_gamma}
\end{minipage}
\begin{minipage}{0.3\textwidth}
\centering
\begin{tikzpicture}
\filldraw[draw=black, thick, fill=gray!20!white] (0,3) -- (0,0) arc(180:0:.5) -- (1,0)arc(180:0:.5) -- (2,3);
\path (0,-.3) node{$0$};
\path (1,-.3) node{$1$};
\path (2,-.3) node{$2$};
\path (-0.3,1) node{$a$};
\path (2.3,1) node{$a$};
\path (0.5,.2) node{$b$};
\path (1.5,.2) node{$b$};
\draw (0,1) arc(90:60:1);
\draw (.5,.866) -- (.5,.5);
\draw (.5,.866) arc(120:60:1);
\draw (1.5,.866) -- (1.5,.5);
\draw (1.5,.866) arc(120:90:1);
\filldraw[fill=black] (0,1) circle[radius=2pt];
\filldraw[fill=black] (1,1) circle[radius=2pt];
\filldraw[fill=black] (2,1) circle[radius=2pt];
\filldraw[fill=black] (.5,.5) circle[radius=2pt];
\filldraw[fill=black] (1.5,.5) circle[radius=2pt];
\filldraw[fill=white] (.5,.866) circle[radius=2pt];
\filldraw[fill=white] (1.5,.866) circle[radius=2pt];
\end{tikzpicture}
\subcaption{$G=\Gamma(2)$.}\label{fig:p_for_gamma_2}
\end{minipage}
\begin{minipage}{0.3\textwidth}
\centering
\begin{tikzpicture}
\filldraw[draw=black, thick, fill=gray!20!white] (0,3) -- (0,0) arc(180:0:.5) -- (1,3);
\draw (0,1) arc(90:60:1);
\draw (.5,.866) -- (.5,.5);
\draw (.5,.866) arc(120:90:1);
\filldraw[fill=black] (0,1) circle[radius=2pt];
\filldraw[fill=black] (1,1) circle[radius=2pt];
\filldraw[fill=black] (.5,.5) circle[radius=2pt];
\filldraw[fill=white] (.5,.866) circle[radius=2pt];
\path (0,-.3) node{$0$};
\path (1,-.3) node{$1$};
\path (-0.3,1) node{$a$};
\path (1.3,1) node{$a$};
\end{tikzpicture}
\subcaption{$G=\Gamma_0(2)$.}\label{fig:p_for_gamma_bot_0_2}
\end{minipage}
\begin{minipage}{0.3\textwidth}
\centering
\begin{tikzpicture}
\filldraw[draw=black, thick, fill=gray!20!white] (0,3) -- (0,0) arc(180:0:.5) -- (1,3);
\draw (0,1) arc(90:60:1);
\draw (.5,.866) -- (.5,.5);
\draw (.5,.866) arc(120:90:1);
\filldraw[fill=black] (0,1) circle[radius=2pt];
\filldraw[fill=black] (1,1) circle[radius=2pt];
\filldraw[fill=black] (.5,.5) circle[radius=2pt];
\filldraw[fill=white] (.5,.866) circle[radius=2pt];
\path (0,-.3) node{$0$};
\path (1,-.3) node{$1$};
\path (-0.3,1) node{$a$};
\path (0.5,.2) node{$a$};
\end{tikzpicture}
\subcaption{$G=\Gamma^0(2)$.}\label{fig:p_for_gamma_top_0_2}
\end{minipage}
\begin{minipage}{0.3\textwidth}
\centering
\begin{tikzpicture}
\filldraw[draw=black, thick, fill=gray!20!white] (-.5,3) -- (-.5,.866) arc(60:0:1) -- (0,0) arc(180:120:1) -- (.5,3);
\draw (-.5,.855) arc(120:60:1);
\filldraw[fill=black] (0,1) circle[radius=2pt];
\filldraw[fill=white] (.5,.866) circle[radius=2pt];
\filldraw[fill=white] (-.5,.866) circle[radius=2pt];
\path (0,-.3) node{$0$};
\end{tikzpicture}
\subcaption{$G=\ker(\Gamma\to\Z/2)$.}\label{fig:p_for_ker_to_z2}
\end{minipage}
\begin{minipage}{0.3\textwidth}
\centering
\begin{tikzpicture}
\filldraw[draw=black, thick, fill=gray!20!white] (0,3) -- (0,0) arc(180:0:.5) -- (1,3);
\draw (0,1) arc(90:60:1);
\draw (.5,.866) -- (.5,.5);
\draw (.5,.866) arc(120:90:1);
\filldraw[fill=black] (0,1) circle[radius=2pt];
\filldraw[fill=black] (1,1) circle[radius=2pt];
\filldraw[fill=black] (.5,.5) circle[radius=2pt];
\filldraw[fill=white] (.5,.866) circle[radius=2pt];
\path (0,-.3) node{$0$};
\path (1,-.3) node{$1$};
\end{tikzpicture}
\subcaption{$G=\ker(\Gamma\to\Z/3)$.}\label{fig:p_for_ker_to_z3}
\end{minipage}
\caption{Special polygons for the subgroups $G\subset\Gamma$ from example~\ref{example_graphs}.}\label{fig:example_graphs_contd}
\end{figure}
\end{example}

\begin{example}
The next two examples, shown in figure~\ref{fig:more_special_poligs}, are more complicated and were calculated using the method described in section~\ref{classical}. %The distinguished edge $[\ii,\rr]$ is shown in bold.
\begin{figure}[h!]
\centering
\begin{minipage}{0.495\textwidth}
\centering
\begin{tikzpicture}[scale=1.5]
\filldraw[draw=black, thick, fill=gray!20!white] (0,4.5) -- (0,0) arc(180:0:1) arc(180:0:1/3) arc(180:0:1/6) arc(180:0:1/2) -- (4,4.5);
\draw (0,4) arc(90:60:4);
\draw (2,3.4641) arc(120:90:4);
\draw (2,3.4641) -- (2,2);
\draw (2,2) -- (2,1.1547);
\draw (2,1.1547) arc to [small, counter clockwise,radius=1.3333] (1.6,.8);
\draw (2,1.1547) arc to [small, clockwise,radius=1.3333] (2.4,.8);
\draw (2.4,.8) arc to [small, clockwise,radius=1.3333] (2.5714,0.4949);
\draw (2.5714,0.4949) arc to [small, counter clockwise,radius=.8] (2.4615,0.3077);
\draw (2.5714,0.4949) arc to [small, clockwise,radius=.5] (2.8,0.4);
\draw (2.8,0.4) arc to [small,clockwise,radius=.5] (2.9231,0.2665);
\draw (2.9231,0.2665) arc to [small,counter clockwise,radius=0.5714] (2.88,.16);
\draw (2.9231,0.2665) arc to [small,clockwise,radius=0.2667] (3.0588,0.2353);
\path (0,-.3) node{$0$};
\path (2,-.3) node{$\frac{1}{2}$};
\path (8/3,-.3) node{$\frac{2}{3}$};
\path (3,-.3) node{$\frac{3}{4}$};
\path (4,-.3) node{$1$};
\path (-0.3,4) node{$a$};
\path (4.3,4) node{$a$};
\path (1.6-.1,.8-.1) node{$b$};
\path (2.88-.05,.16-.15) node{$b$};
\path (2.4615-.05,0.3077-.14) node{$c$};
\path (3.0588+.14,0.2353) node{$c$};
\filldraw[fill=black] (2.88,.16) circle[radius=1pt];
\filldraw[fill=black] (3.0588,0.2353) circle[radius=1pt];
\filldraw[fill=black] (2.8,0.4) circle[radius=1pt];
\filldraw[fill=black] (2.4615,0.3077) circle[radius=1pt];
\filldraw[fill=black] (0,4) circle[radius=1pt];
\filldraw[fill=black] (4,4) circle[radius=1pt];
\filldraw[fill=black] (2,2) circle[radius=1pt];
\filldraw[fill=black] (1.6,.8) circle[radius=1pt];
\filldraw[fill=black] (2.4,.8) circle[radius=1pt];
\filldraw[fill=white] (2,3.4641) circle[radius=1pt];
\filldraw[fill=white] (2,1.1547) circle[radius=1pt];
\filldraw[fill=white] (2.5714,0.4949) circle[radius=1pt];
\filldraw[fill=white] (2.9231,0.2665) circle[radius=1pt];
\end{tikzpicture}
\subcaption{$G=\Gamma_0(6)$.}\label{fig:spec_poly_gamma_6}
\end{minipage}
%\hspace{2cm}
\begin{minipage}{0.495\textwidth}
\centering
\begin{tikzpicture}[scale=1.5]
\filldraw[draw=black, thick, fill=gray!20!white] (0,4.5) -- (0,0) arc to [small,clockwise,radius=.5714] (1.0769,.2665) arc to [small,clockwise,radius=.2667] (4/3,0) arc to [small,clockwise,radius=1/3] (2,0) arc to [small,clockwise,radius=1/3] (8/3,0) arc to [small,clockwise,radius=.2667]  (2.9230,0.2665) arc to [small,clockwise,radius=.5714] (4,0) -- (4,4.5);
\draw (0,4) arc(90:60:4);
\draw (2,3.4641) arc(120:90:4);
\draw (2,3.4641) -- (2,2);
\draw (2,2) -- (2,1.1547);
\draw (2,1.1547) arc to [small, counter clockwise,radius=1.3333] (1.6,.8);
\draw (2,1.1547) arc to [small, clockwise,radius=1.3333] (2.4,.8);
\draw(1.6,.8) arc to [small,counter clockwise, radius=1.3333] (1.4286,0.4949);
\draw (2.4,.8) arc to [small, clockwise,radius=1.3333] (2.5714,0.4949);
\draw (1.4286,0.4949) arc to [small, clockwise,radius=.8] (1.5385,0.3077);
\draw (2.5714,0.4949) arc to [small, counter clockwise,radius=.8] (2.4615,0.3077);
\draw (1.4286,0.4949) arc to [small,counter clockwise, radius=.5] (1.2,.4);
\draw (2.5714,0.4949) arc to [small, clockwise,radius=.5] (2.8,0.4);
\draw (1.2,.4) arc to [small, counter clockwise, radius=.5] (1.0769,0.2665);
\draw (2.8,0.4) arc to [small, clockwise, radius=.5] (2.9231,0.2665);
\filldraw[fill=white] (1.4286,0.4949) circle[radius=1pt];
\filldraw[fill=white] (2.5714,0.4949) circle[radius=1pt];
\filldraw[fill=white] (2,1.1547) circle[radius=1pt];
\filldraw[fill=white] (2,3.4641) circle[radius=1pt];
\filldraw[fill=white] (1.0769,0.2665) circle[radius=1pt];
\filldraw[fill=white] (2.9230,0.2665) circle[radius=1pt];
\filldraw[fill=black] (1.5385,0.3077) circle[radius=1pt];
\filldraw[fill=black] (2.4615,0.3077) circle[radius=1pt];
\filldraw[fill=black] (0,4) circle[radius=1pt];
\filldraw[fill=black] (4,4) circle[radius=1pt];
\filldraw[fill=black] (2,2) circle[radius=1pt];
\filldraw[fill=black] (1.6,.8) circle[radius=1pt];
\filldraw[fill=black] (2.4,.8) circle[radius=1pt];
\filldraw[fill=black] (1.2,.4) circle[radius=1pt];
\filldraw[fill=black] (2.8,0.4) circle[radius=1pt];
\path (0,-.3) node{$0$};
\path (4/3,-.3) node{$\frac{1}{3}$};
\path (2,-.3) node{$\frac{1}{2}$};
\path (8/3,-.3) node{$\frac{2}{3}$};
\path (4,-.3) node{$1$};
\path (-0.3,4) node{$a$};
\path (4.3,4) node{$a$};
\end{tikzpicture}
\subcaption{$G=\Gamma_0(13)$.}\label{fig:spec_poly_gamma_13}
\end{minipage}
\caption{More examples of special polygons.}\label{fig:more_special_poligs}
\end{figure}
\end{example}

\begin{definition}
A subset $\{g_1,\ldots,g_l\}$ of a group $H$ is called an {\it independent set of generators} if $H$ is the free product of the cyclic subgroups generated by the $g_i$'s. 
\end{definition}

%Now we briefly describe how to obtain an independent set of generators for $G\subset\Gamma$ from the above information. 
\begin{theorem}\label{main1}
\begin{enumerate}
\item The couple $(P,\sigma)$ is a special polygon.
\item Let $S$ be a subset of the set of the edges of $P$ which intersects every $\sigma$-orbit exactly once. Then $\{g_s\mid s\in S\}$ is an independent set of generators of $G$.
\end{enumerate}
\end{theorem}

Observe that in the theorem we do not assume that the index $(\Gamma:G)<\infty$. Similarly, the proof of Theorem on p.~1055 in~\cite{kulkarni} extends verbatim to arbitrary subgroups of $\Gamma$.

{\bf Proof.} The first assertion follows from the construction. The second one follows from~\cite[I.5, theorem 13]{serre} or~\cite[theorem on p.~1055]{kulkarni}.$\clubsuit$

\bigskip

R.~Kulkarni noticed~\cite[Theorem 4.2]{kulkarni} that the isomorphism classes of bipartite cuboid graphs correspond bijectively to the conjugacy classes of subgroups of $\Gamma$. (The theorem is stated only for finite index subgroups but is true without this assumption.) The mapping that associates a graph to a conjugacy class of subgroups is constructed in \cite{kulkarni} as follows: we take a conjugacy class $C$ to the isomorphism class of the image of $\mathcal{T}$ in $G\backslash\mathbb{H}$ where $G\in C$ is a representative. It follows that $\Xc(G)$ is a representative in the isomorphism class of bipartite cuboid graphs that corresponds to the conjugacy class of $G$.

A natural question is how to describe the set of subgroups of $\Gamma$, rather than the set of their conjugacy classes, in terms of graphs. 

\begin{definition}
A {\it pointed} bipartite cuboid graph is the data of a bipartite cuboid graph and an edge of the graph.
\end{definition}

\begin{theorem}\label{pointed}
The map $G\mapsto(\Xc(G),\mathbf{r}_G)$ induces a bijection between the set of subgroups of $\Gamma$ and the set of isomorphism classes of pointed bipartite cuboid graphs. In particular, the quotient of the set $\mathop{\mathrm{Edge}}(\Xc(G))$ by the group $\mathop{\mathrm{Aut}}\Xc(G)$ of the automorphisms of $\Xc(G)$ that preserve the bipartite cuboid structure is in bijection with the conjugacy class of $G$; so $G$ is normal in $\Gamma$ iff $\mathop{\mathrm{Aut}}\Xc(G)$ acts transitively on $\mathop{\mathrm{Edge}}(\Xc(G))$.
\end{theorem}

{\bf Proof.} Suppose that $G_1, G_2$ are two subgroups such that the graphs $\Xc(G_1), \Xc(G_2)$ are isomorphic via an isomorphism that takes the distinguished edge to the distinguished edge. Then we may assume that the trees $Y'_1, Y_2'$ obtained by cutting $\Xc(G_1)$, respectively $\Xc(G_2)$ along some of the vertices of type~(0) are the same. Let us embed $Y'_1=Y'_2$ in $\mathbf{T}$ so that the distinguished edge goes to the distinguished edge of $\mathbf{T}$. Then the diagram~(\ref{diag}) commutes, see remark~\ref{rmk:what_distinguished_edges_are_needed_for}. 
%Hindsight
By constructing the corresponding special polygon as described earlier in this section and applying part 2 of theorem~\ref{main1} we deduce that the subgroups $G_1$ and $G_2$ are generated by the same subset of $\Gamma$.$\clubsuit$

\begin{example}
Observe that for the graphs shown in figure~\ref{fig:1}, the group $\mathop{\mathrm{Aut}}\Xc(G)$ is transitive on $\mathop{\mathrm{Edge}}(\Xc(G))$ in all cases except for $G=\Gamma_0(2)=\Gamma_1(2)$ and $G=\Gamma^0(2)=\Gamma^1(2)$.
\end{example}

\begin{remark}\label{red}
Let $(P,\sigma)$ be a special polygon for a finite index subgroup $G\subset\Gamma$. One can then find for a given $z\in\Hh$ a point $w\in P$ and an element $g\in G$ such that $gz=w$. Let us sketch the algorithm. Let $z_0$ be an interior point of $P$, and set $l$ to be the geodesic \emph{line} containing $\gamma=[z_0,z]$. Let us orient $\gamma$ from $z_0$ to $z$.

The tessellation of $\Hh$ by the $G$-copies of $P$ subdivides $\gamma$ into geodesic line segments $\gamma_0,\ldots,\gamma_m$, which we assume ordered according to the orientation of $\gamma$. Let us construct these line segments, and also elements $g_i\in G,i=0,\ldots,m$ such that $g_i\gamma_i\subset P$ for all $i$. We have $\gamma_0=\gamma\cap P$, and we set $g_0=\Id$. If we know $\gamma_i$ and $g_i$, we can find an element $h\in G$ such that $g_{i+1}=hg_i$ takes $\gamma_{i+1}$ to $P$, and either $h$ or $h^{-1}$ belongs to the system of independent generators $\{g_s\mid s\in S\}$ from part~2 of theorem~\ref{main1}. This can be done as follows. Let $z'\in \partial P$ be the end point of $g_i\gamma_i$. There are three cases to consider:

\begin{itemize}
\item $z'$ is not a vertex of $P$. Then we set $h\in G$ to be the element such that $h$ takes the side $s\ni z'$ to the side $\sigma(s)$ of $P$.
\item $z'$ is an elliptic vertex of $P$ of order 2. Then we set $h\in G$ to be the element of order 2 such that $h$ fixes $z'$ and interchanges the sides of $P$ that meet at $z'$.
\item $z'$ is an elliptic vertex of $P$ of order 3. Let $s_1,s_2$ be the sides of $P$ meeting at $z'$. The angle $\alpha$ between $g_i\gamma_i$ and one of these sides, say $s_1$, is $\leq\frac{\pi}{3}$. We set $h\in G$ to be the element of order 3 such that $h$ fixes $z'$ and takes $s_2$ to $s_1$. (It may happen that $\alpha=\frac{\pi}{3}$. In this case $g_{i+1}\gamma_{i+1}$ will be inside $s_2$, and the procedure terminates at $i+1$, i.e.\ we have $m=i+1$.)
\end{itemize}

Once we have constructed $g_{i+1}$ we can recover $\gamma_{i+1}$ as $g_{i+1}^{-1}((g_{i+1}l)\cap P)$. After the procedure terminates we set $g=g_m$ and $w=gz$.

Note that this algorithm returns $g$ as a product of the independent generators $\{g_s\mid s\in S\}$ and their inverses. So we can write a given $g\in G$ in terms of these generators by applying the algorithm to $z=g z_0$.
\end{remark}

\section{Automorphism groups}\label{aut_gr}

Each subgroup $G\subset\Gamma$ gives rise to several automorphism groups. These are the automorphism groups $\Aut(Y(G))$ and $\Aut(X(G))$ of the smooth complex curves $Y(G)=G\backslash\Hh$ and $X(G)=G\backslash \Hh^*$ (see e.g.\ \cite[section 2.4]{diashur}); the isometry group $\Iso(Y(G))$ of $Y(G)$ viewed as a complete hyperbolic surface\footnote{This means the corresponding metric space is complete.}, possibly with conical singularities; the automorphism group $\Aut(\Xc(G))$ of the bipartite cuboid graph $\Xc(G)$; and the groups $N_{\PSL_2(\R)}(G)/G$ and $N_\Gamma(G)/G$ where $N_H(G)$ denotes the normaliser of $G$ in $H$. Here we discuss the relationship between these groups, and the algorithms for calculating them that we are aware of. Let us start with the graph $\Xc(G)$.

\begin{lemma} 
If $H_1$ is a group and $H_2\subset H_1$ is a subgroup, then the automorphism group of $H_2\backslash H_1$ as a right $H_1$-set is $N_{H_1}(H_2)/H_2$.
\end{lemma}

$\clubsuit$

\begin{Prop}\label{aut_graph}
We have $\Aut(\Xc(G))=N_\Gamma(G)/G$.
\end{Prop}

{\bf Proof.} An automorphism $a$ of $\Xc(G)$ induces a bijection $a_{\mathrm{Edge}}$ of the set of the edges, which is $G\backslash\Gamma$, and a compatible bijection of the set of the vertices $(G\backslash\Gamma/G_0)\sqcup (G\backslash\Gamma/G_1)$ that preserves both summands of the disjoint union. The map $a_{\mathrm{Edge}}:G\backslash\Gamma\to G\backslash\Gamma$ takes two element $G_0$-orbits to two element $G_0$-orbits, and similarly for one element $G_0$-orbits. So $a$ commutes with the right action of $G_0$. A similar argument combined with the fact that $a$ preserves the cyclic order at each vertex shows that $a_{\mathrm{Edge}}$ commutes with the right action of $G_1$. So $a_{\mathrm{Edge}}$ is an automorphism of $G\backslash\Gamma$ as a right $\Gamma$-set. Applying the lemma we conclude that $a\in N_{\Gamma}(G)/G$.$\clubsuit$

\medskip

We now turn to the complex curve $Y(G)$. We do not know how to calculate the full automorphism group $\Aut(Y(G))$. But here is some partial information. Let $\Aut^{\mathrm{ell}}(Y(G))\subset\Aut(Y(G))$ be the subgroup of the automorphisms of $Y(G)$ that preserve the elliptic points and their types. This subgroup need not coincide with $\Aut(Y(G))$, as can be seen by taking $G=\Gamma$, although it clearly does when $Y(G)$ has no elliptic points, i.e.\ when $G$ is torsion free. Let $\Iso^+(Y(G))\subset \Iso(Y(G))$ denote the group of isometries that preserve orientation, i.e.\ act trivially on $H^2_c(Y(G),\Z)$.

\begin{Prop}\label{iso_aut_ngg}
We have $\Iso^+(Y(G))=\Aut^{\mathrm{ell}}(Y(G))=N_{\PSL_2(\R)}(G)/G$.
\end{Prop}

{\bf Proof.} %The inclusion $\Iso^+(Y(G))\subset \Aut^{\mathrm{ell}}(Y(G))$ follows from the Myers-Steenrod theorem and basic hyperbolic geometry.
Let us first show that $\Iso^+(Y(G))\subset \Aut^{\mathrm{ell}}(Y(G))$. By the Myers-Steenrod theorem, an isometry $f$ of $Y(G)$ as a metric space is smooth on $Y'(G)=Y(G)$ minus the elliptic points, and preserves the hyperbolic metric on $Y'(G)$. So if $f$ preserves orientation, it induces a holomorphic automorphism of $Y'(G)$.

To handle the elliptic points if any we use these observations: An elliptic point of order $l\in\{2,3\}$ has an open neighbourhood isometric to $V_\varepsilon=H\backslash U_\varepsilon$ where $U_\varepsilon\subset\C$ is the disk of radius $\varepsilon>0$ about the origin with the Poincar\'e disk metric, and $H\subset\C^*$ is the group of $l$-th roots of $1$. The orientation preserving isometries of $V_\varepsilon$ are induced by the rotations about the origin, and the holomorphic coordinate in $V_\varepsilon$ is induced by the function $z\mapsto z^l$ on $U_\varepsilon$, which makes the isometries of $V_\varepsilon$ holomorphic.

Next we will show that $\Aut^{\mathrm{ell}}(Y(G))\subset N_{\PSL_2(\R)}(G)/G$. Let $p:\Hh\to Y(G)$ be the projection. Suppose $f:Y(G)\to Y(G)$ is an automorphism that preserves the set $E\subset Y(G)$ of elliptic points, and the type of each element of $E$. 
%Choose a base point $y_0\in Y(G)$. The subgroup of $\pi_1(Y(G), y_0)$ that corresponds to the covering $p|_{\Hh\backslash p^{-1}(E)}:\Hh\backslash p^{-1}(E)\to Y(G)\backslash E$ is the normal closure of the set of loops that start at $y_0$, go along a path $\gamma$ to a point $e'$ close to an elliptic point $e\in E$, make two or three turns around $e$ and then return to $y_0$ along $\gamma^{-1}$. 
By the standard theory of covering spaces, $f|_{Y(G)\backslash E}$ lifts to a holomorphic automorphism $\tilde f:\Hh\backslash p^{-1}(E)\to \Hh\backslash p^{-1}(E)$. 
%Moreover, the preimage under $p$ of a sufficiently small neighbourhood $U\cong \{z\in\C||z|<1\}$ of a point $e\in E$ is analytically isomorphic to a disjoint union of copies of $U$. 
Moreover, the singularity of $\tilde f$ at each point of $p^{-1}(E)$ is removable, so $\tilde f$ extends to a holomorphic automorphism $\tilde f:\Hh\to\Hh$ that covers $f$. 
%indeed, post-compose $\tilde f$ with an isomorphism $\mathbb{H}\to $ open unit disk. The composition is globally bounded, so all singularities are removable. Hence the same is true for $\tilde f$.
For every $x\in\Hh, g\in G$ there is an $A(x,g)\in G$ such that $$\tilde f(gx)=A(x,g)\tilde f(x).$$ Since $G$ is discrete, $A(x,g)$ depends only on $g$ and not on $x$. We conclude that $\tilde f g {\tilde f}^{-1}\in G$ for all $g\in G$, which implies $f\in N_{\PSL_2(\R)}(G)/G$.

Finally, an element of $N_{\PSL_2(\R)}(G)/G$ clearly preserves both the metric on $Y(G)$ and the orientation induced from $\Hh$, showing that $N_{\PSL_2(\R)}(G)/G\subset \Iso^+(Y(G))$.$\clubsuit$

\smallskip

Of course, in general $N_\Gamma(G)/G\neq N_{\PSL_2(\R)}(G)/G$, but here is a sufficient condition for these two groups to coincide. For a cusp $x\in X(G)$, the edges of the copies of the triangle $\triangle$ cut a small neighbourhood of $x$ into some number of sectors. We will call this number the {\it width} of the cusp $x$, denoted $w_x$. Each copy of $\triangle$ contributes 0, 1, or 2 to $w_x$.

\begin{Prop}\label{equal_width}
If the widths of all cusps of $X(G)$ are the same, then $N_\Gamma(G)/G= N_{\PSL_2(\R)}(G)/G$.
\end{Prop}

For example, this proposition is applicable to $G=\Gamma(N)$ as the group $\SL_2(\Z/N)/\{\pm I\}$ acts transitively on the cusps of $X(N)$, which gives an alternative proof of~\cite[Corollary 3.7]{bkx}.

{\bf Proof.} In the proof we will view $Y(G)$ as a complete hyperbolic surface, with conical singularities at the elliptic points. For a cusp $x$ of $X(G)$, let $\tilde x$ be a lift of $x$ to the closure $\bar P$ (in $\bar\Hh=\Hh\sqcup\R\sqcup\{\infty\}$) of the fundamental domain $P$ from section~\ref{grfun}. Let $s_1,s_2$ be the geodesic lines that contain the sides of $P$ meeting at $\tilde x$. For $t>0$ set $A_{t,\tilde x}$ to be the closed subset of $\Hh$ of area $\frac{t}{w_x}$ bounded by $s_1, s_2$ and a horocycle centred at $\tilde x$. Let $A_{t,x}=\bigcup_{\tilde x} A_{t,\tilde x}$ where the union is taken over all lifts $\tilde x$ of $x$ to $\bar P$.  The image $B_{t,x}$ of $A_{t,x}$ in $Y(G)$ will be called the {\it closed horocyclic neighbourhood} of $x$ of area $t$. 

Let $B_t$ be the union of $B_{t,x}$ for all cusps $x$ of $X(G)$. Observe that every $B_t$ is stable under all automorphisms of $Y(G)$. Let us see how one can reconstruct the graph $\Xc(G)\subset Y(G)$ from $B_t$. Let $t_0$ be the smallest value of $t$ for which $B_{t,x}$ start to (self-)intersect, and let $t_1$ be the smallest value of $t$ such that $B_t=Y(G)$. (One can calculate $t_0=w\int_{D_0}\frac{dx\wedge dy}{y^2}=\frac{w}{2}, t_1=w\int_{D_1}\frac{dx\wedge dy}{y^2}=\frac{w}{\sqrt{3}}$ where $D_0$ and $D_1$ are given by $0\leq x\leq\frac{1}{2}, y\geq 1$, respectively by $0\leq x\leq\frac{1}{2}, y\geq \frac{\sqrt{3}}{2}$, and $w$ is the common width of all cusps of $X(G)$.)

Let $\partial B_t$ denote the topological boundary of $B_t$, i.e.\ the closure minus the interior. For $t\in [t_0,t_1)$ let $V(t)$ be the set of non-smooth points of $\partial B_t$\footnote{These are the points of $\partial B_t$ at which $Y(G)$ has a singularity, and the smooth points $y\in Y(G)$ which belong to $\partial B_t$ and such that locally at $y$ the set $\partial B_t$ is not the zero locus of a smooth function with non-vanishing differential.}, and let $V(t_1)$ be the set $\bar V\backslash V$ where $V=\bigcup_{t_0\leq t< t_1} V(t)$. Then $\Xc(G)=\bigcup_{t_0\leq t\leq t_1} V(t)$, and $V(t_0)$, respectively $V(t_1)$ is the set of all vertices of $\Xc(G)$ of type~(0), respectively of type~(1). To see this let $\bar\triangle$ be the closure of $\triangle$ in $\bar\Hh$, and recall that by proposition~\ref{fund_domain_from_graph} the polygon $\bar P$ is the union of the images of $\bar\triangle$ by elements of $\Gamma$. If we take any such copy $\bar\triangle\to \bar P$ and take the preimage of $B_t$ under the composite $\bar\triangle\to\bar P\to X(G)$, we get $\bar\triangle\cap(\bar H_0\cup \bar H_\infty)$ where $H_0$ and $H_\infty$ are horodisks centred at $0$, respectively $\infty$, and such that $\triangle\cap H_0$ and $\triangle\cap H_\infty$ have equal area. 

We conclude that every automorphism $f$ of the hyperbolic surface $Y(G)$ preserves $\Xc(G)$, which implies by proposition~\ref{aut_graph} that $f\in N_\Gamma(G)$.$\clubsuit$

\medskip
Finally, we will discuss the compact complex curve $X(G)$. The group $\Aut(X(G))$ is more mysterious than those considered so far, and it seems that only partial results about $\Aut(X(G))$ are known. The automorphism group $\Aut(Y(G))$ of the non-compact smooth curve $Y(G)$ is included in $\Aut(X(G))$ as the subgroup that preserves the set of the cusps. An automorphism of $X(G)$ that does not belong to $N_{\PSL_2(\R)}(G)/G$ is called {\it exceptional}. These exist but seem to be rare. Below we briefly review the results for the classical congruence subgroups.

For $G=\Gamma_0(N)$, the group $N_{\PSL_2(\R)}(G)$ is given in several sources, see e.g.\ J.~Lehner and M.~Newman~\cite{ln}, or J.~H.~Conway and S.~P.~Norton~\cite{cn}. The quotient group $N_{\PSL_2(\R)}(G)/G$ was calculated by A.~O.~L.~Atkin and J.~Lehner~\cite[\S 4]{al} with corrections by M.~Akbas and D.~Singerman~\cite{as}, see also F.~Bars~\cite{bars}. M.~A.~Kenku and F.~Momose proved~\cite{ken_mom} with complements and corrections by N.~D.~Elkies~\cite{elkies} and M.~Harrison~\cite{harri} that if the genus of $X_0(N)=\Gamma_0(N)\backslash\Hh^*$ is $>1$, then $X_0(N)$ has an exceptional automorphism if and only if $N=37, 63$, or $108$. 

For $G=\Gamma_1(N)$, the group $N_{\PSL_2(\R)}(G)$ was determined by M.-L.\ Lang~\cite{lang_norma}. The group $N_{\PSL_2(\R)}(G)/G$ in this case is an extension of $(\Z/2)^{|\{\mbox{{\scriptsize prime divisors of }}N\}|}$ by $(\Z/N)^*$, see S.~Zemel~\cite[Corollary 3.3]{zemel}, but a full description of this group appears to be unknown. It seems to be folklore knowledge that $X_1(N)=\Gamma_1(N)\backslash\Hh^*$ does not have exceptional automorphisms if the genus is $>1$, but we were unable to locate a reference for this. 

Finally, for $G=\Gamma(N)$ all automorphism groups coincide:

\begin{Prop}\label{aut_gamma_n}
For $G=\Gamma(N), N\geq 7$ we have
$$\Aut (\Xc(G))=\Iso^+(Y(G))=\Aut(Y(G))=\Aut (X(G))=\SL_2(\Z/N)/\{\pm I\}.$$
\end{Prop}

This is Theorem~3.1 in~\cite{bkx}, although the result may have been known earlier, see the discussion in the introduction of that paper. Below we give an alternative proof, which we believe is simpler.

{\bf Proof.} The first two equalities follow from propositions~\ref{aut_graph}-\ref{equal_width}. Further, note that the group $\SL_2(\Z/N)/\{\pm I\}$ acts on $X(N)=\Gamma(N)\backslash\Hh^*$ by automorphisms, and this action preserves $Y(N)=\Gamma(N)\backslash\Hh$. So to prove the proposition it suffices to show that any automorphism of $X(N)$ is induced by an element of $\SL_2(\Z/N)/\{\pm I\}$.

Suppose that a finite group $H$ acts on a compact 2-dimensional smooth manifold $S$ by diffeomorphisms, and assume the action is cellular with respect to some CW-structure on $S$ and that the point-wise and set-wise stabilisers of each cell are the same. One can then define (see~W.~Thurston~\cite[Chapter 13, Definition 13.3.3]{thurston}) the {\it orbifold Euler characteristic} $\chi^{orb}(S/H)$ of $S/H$ by counting every cell $e$ with coefficient $(-1)^{\dim e} \frac{1}{|\Stab(e)|}$ where $\Stab(e)$ is the point-wise stabiliser of $e$. Assume for simplicity that there are finitely many points with non-trivial stabiliser, and that all stabilisers are orientation preserving, hence cyclic. We have then
\begin{equation}\label{orbi_euler}
\chi^{orb}(S/H)=\chi(S/H)-\sum_i \frac{n_i-1}{n_i}
\end{equation}
where the sum is over all non-free orbits and $n_i$ is the cardinality of the stabiliser of a point of the $i$-th orbit (ibid., formula 13.3.4). The number $\chi^{orb}(S/H)$ is not necessarily an integer, but we have
\begin{equation}\label{orbi_quotient}
\chi^{orb}(S/H)=\frac{\chi(S)}{|H|}
\end{equation}
regardless of whether or not the action is free (ibid., Proposition 13.3.4).

\smallskip

The quotient of $X(N)$ by $\SL_2(\Z/N)/\{\pm I\}$ is $X(1)=\Gamma\backslash\Hh^*\cong\mathbb{P}^1(\C)$. There are three types of non-free orbits for the action of $\SL_2(\Z/N)/\{\pm I\}$ on $X(N)$, namely the images of $\ii$ and $\rr$ in $X(1)$, and the cusp of $X(1)$: over all other points of $X(1)$ the covering $X(N)\to X(1)$ is unramified. The corresponding stabilisers are $\Z/2, \Z/3$, and $\Z/N$ respectively. The first two of these assertions follow from the fact that the action of $\SL_2(\Z/N)/\{\pm I\}$ comes from $N_\Gamma(\Gamma(N))/\Gamma(N)$, so it permutes the images in $X(N)$ of the copies of the closure $\bar \triangle\subset\Hh^*$. To see the third one, observe that the stabiliser in $\SL_2(\Z/N)/\{\pm I\}$ of the cusp at infinity is generated by the image of $\left[\begin{smallmatrix}1&1\\0&1\end{smallmatrix}\right]\in\Gamma$.

Now using~(\ref{orbi_euler}) we conclude that $$\chi^{orb}\Bigl(X(N)/(\SL_2(\Z/N)/\{\pm I\})\Bigr)=2-\left(1-\frac{1}{2}\right)-\left(1-\frac{1}{3}\right)-\left(1-\frac{1}{N}\right)=-\frac{1}{6}+\frac{1}{N}.$$ %(Observe that formula~(\ref{orbi_quotient}) gives the same answer.) 
A straightforward combinatorial argument (see lemma~\ref{lemma_comb_check}) shows that there are finitely many collections $N, n_1,\ldots, n_k$ of integers such that $N\geq 7$, all $n_i\geq 2$, and $$\frac{-\frac{1}{6}+\frac{1}{N}}{2-\sum_{i=1}^k \frac{n_i-1}{n_i}}$$ is an integer $\geq 2$. For each of these, none of the $n_i$ is divisible by $N$. On the other hand, if we have a collection that comes from an action of a finite group containing $\SL_2(\Z/N)/\{\pm I\}$ on $X(N)$ by automorphisms, then one of the stabilisers contains $\Z/N$, and $-\frac{1}{6}+\frac{1}{N}$ is a positive integral multiple of $\chi^{orb}$ of the quotient by formula~(\ref{orbi_quotient}). We conclude that $\SL_2(\Z/N)/\{\pm I\}$ cannot be contained in a strictly larger finite automorphism group of $X(N)$. Since we assume $N\geq 7$, the genus of $X(N)$ is at least $2$, so $\Aut(X(N))$ is finite. This implies $\Aut(X(N))=\SL_2(\Z/N)/\{\pm I\}$.$\clubsuit$

\begin{remark} 
For $N\leq 6$ the genus of $X(N)$ is $\leq 1$, so $\Aut(X(N))$ is infinite. The curve $X(7)$ is the unique (up to isomorphism) curve of genus $3$ with $168=84(3-1)$ automorphisms, which is the maximal number allowed by the Hurwitz bound. So $X(7)$ is in fact isomorphic to the Klein quartic.
\end{remark}

\begin{lemma}\label{lemma_comb_check}
There are finitely many collections $N, n_1,\ldots, n_k$ of integers such that $N\geq 7$, $2\leq n_1\leq\cdots \leq n_k$, and 
\begin{equation}\label{eq_d}
d=\frac{-\frac{1}{6}+\frac{1}{N}}{2-\sum_{i=1}^k \frac{n_i-1}{n_i}}
\end{equation}
is an integer $\geq 2$. These collections exist only for $k=3$ and are given in table~(\ref{table:tab_triangle_groups}).
\begin{table}[h!]
\centering
{\tabulinesep=1mm
\begin{tabu}{|c|c|c|}
\hline
$(n_1,n_2,n_3)$ & $N$ & $d$\\
\hline
$(2,3,7)$ & $14$ & $4$\\
\hline
$(2,3,7)$ & $21$ & $5$\\
\hline
$(2,3,7)$ & $42$ & $6$\\
\hline
$(2,3,8)$ & $12$ & $2$\\
\hline
$(2,3,8)$ & $24$ & $3$\\
\hline
$(2,3,9)$ & $18$ & $2$\\
\hline
$(2,3,10)$ & $30$ & $2$\\
\hline
$(2,3,11)$ & $66$ & $2$\\
\hline
$(2,4,5)$ & $15$ & $2$\\
\hline
$(2,4,5)$ & $60$ & $3$\\
\hline
\end{tabu}
}
\caption{Integer solutions of equation~(\ref{eq_d}) such that $N\geq 7$, $2\leq n_1\leq\cdots \leq n_k$, $d\geq 2$.}
\label{table:tab_triangle_groups}
\end{table}

\end{lemma}

{\bf Proof.} Let $N, n_1,\ldots, n_k$ be a collection that satisfies the conditions of the lemma, and define $d$ using formula~(\ref{eq_d}) and assuming the denominator is non zero. Set
$$S=2-\sum_{i=1}^k \frac{n_i-1}{n_i}.$$ 

{\bf Observation.} If $|S|\geq \frac{1}{12}$ for some choice of integers $2\leq n_1\leq\cdots \leq n_k$, then for these $n_1,\ldots, n_k$ and any $N\geq 7$ we will have $d<2$. 

Now, to begin with, note that $k$ must be $\geq 3$ as otherwise $d<0$. Similarly, if $k\geq 5$, then $|S|\geq |{2-\frac{5}{2}}|=\frac{1}{2}>\frac{1}{12}$.

Next we exclude $k=4$. In this case the minimum non-zero value of $|S|$ is attained for $(n_1,n_2,n_3,n_4)=(2,2,2,3)$ and is equal $\frac{1}{6}>\frac{1}{12}$.

So we must have $k=3$. Suppose $n_1\geq 4$. Again, the minimum value of $|S|$ is attained for $(n_1,n_2,n_3)=(4,4,4)$ and is equal $\frac{1}{4}>\frac{1}{12}$. So $n_1$ must be $2$ or $3$.

\smallskip

{\bf Case $n_1=3$.} We start with the case $n_1=3$. Suppose $n_2\geq 4$. Then the minimum value of $|S|$ is $\frac{1}{6}>\frac{1}{12}$, attained for $(n_1,n_2,n_3)=(3,4,4)$. So $n_2=3$. Similarly, if $n_1=n_2=3$ and $n_3\geq 6$, then $|S|\geq \frac{1}{6}>\frac{1}{12}$. It remains to consider the cases $(n_1,n_2,n_3)=(3,3,3), (3,3,4)$ and $(3,3,5)$. The corresponding values for $|S|$ are $0, \frac{1}{12}$, and $\frac{2}{15}>\frac{1}{12}$ respectively.

\smallskip

{\bf Case $n_1=2$.} Suppose now $n_1=2$. If $n_2=2$, then $S>0$, which makes $d<0$. If $n_2\geq 5$, then the minimum value of $|S|$ is attained for $(n_1,n_2,n_3)=(2,5,5)$ and is equal $\frac{1}{10}>\frac{1}{12}$. So $n_2$ must be $3$ or~$4$.

{\bf Subcase $n_2=4$.} If $n_3=4$, then $S=0$, so $d$ is undefined. If $n_3\geq 6$, then $|S|\geq |2-\frac{1}{2}-\frac{3}{4}-\frac{5}{6}|=\frac{1}{12}$. The remaining possibility is $n_3=5$. In this case $S=-\frac{1}{20}$. There are two integers $d$ such that $2\leq d<\frac{\frac{1}{6}}{\frac{1}{20}}$, namely $d=2$ and $d=3$. Substituting in~(\ref{eq_d}) and solving for $N$ we get $N=15$, respectively $N=60$. This gives us the last two rows of table~(\ref{table:tab_triangle_groups}).

{\bf Subcase $n_2=3$.} Finally, suppose $n_2=3$. If $n_3\leq 6$, then $S\geq 2-\frac{1}{2}-\frac{2}{3}-\frac{5}{6}=0$, which implies that $d<0$ or is undefined. If $n_3\geq 12$, then $|S|\geq \left|2-\frac{1}{2}-\frac{2}{3}-\frac{11}{12}\right|=\frac{1}{12}$. So we need to consider the cases $n_3=7,\ldots, 11$. In each of these we argue as in the case $(n_1,n_2,n_3)=(2,4,5)$ above using the following table.

\begin{center}
{\tabulinesep=1mm
\begin{tabu}{|c|c|c|c|}
\hline
$n_3$ & $S$ & integers $2\leq d< \frac{\frac{1}{6}}{|S|}$ & $N$\\
\hline 
$7$ &$-\frac{1}{42}$&$2$,$3$,$4$,$5$,$6$&$\frac{42}{5}$, $\frac{21}{2}$, $14$, $21$, $42$\\
\hline
$8$ &$-\frac{1}{24}$&$2$,$3$&$12$, $24$\\
\hline
$9$ &$-\frac{1}{18}$&$2$&$18$\\
\hline
$10$ &$-\frac{1}{15}$&$2$&$30$\\
\hline
$11$ &$-\frac{5}{66}$&$2$&$66$\\
\hline
\end{tabu}
}
\end{center}

We then discard the solutions with non-integer $N$ and get the remaining part of table~(\ref{table:tab_triangle_groups}).$\clubsuit$

\bigskip

We will now briefly discuss algorithms for calculating $N_\Gamma(G)$ and $N_{\PSL_2(\R)}(G)$, and the quotients $N_\Gamma(G)/G$ and $N_{\PSL_2(\R)}(G)/G$. So from now on let us assume $(\Gamma:G)<\infty$. In~\cite[section 6]{mong_lung_lang}, M.-L.~Lang gives an algorithm for calculating a system of generators for $N_{\PSL_2(\R)}(G)$ starting from a free system of generators for $G$ (see section~\ref{gen}), and the fundamental polygon (see section~\ref{grfun}). 

More precisely, it is shown in op.\ cit.\ that $N_{\PSL_2(\R)}(G)$ is included in a certain subgroup of $\PSL_2(\R)$, denoted $\Gamma^+(X|Y)$, such that $(\Gamma^+(X|Y):G)<\infty$, and there is an algorithm for giving a complete set $\{y_1,\ldots, y_v\}$ of coset representatives for $\Gamma^+(X|Y)/G$. The normaliser $N_{\PSL_2(\R)}(G)$ is then generated by $G$ and the $y_i$ such that $y_i G y_i^{-1}=G$. Observe that this also gives us a generating set for $N_\Gamma(G)$: we can take the generating set for $N_{\PSL_2(\R)}(G)$ we've just described and remove those $y_i$ that do not belong to~$\Gamma$. 

For understanding the finite groups $N_\Gamma(G)/G$ and $N_{\PSL_2(\R)}(G)/G$ having a faithful linear representation can be more useful than generators and relations.

\begin{lemma}\label{fin_auto_faithful}
Let $X$ be a smooth connected complex algebraic curve of negative Euler characteristic. Then the automorphism group $\Aut X$ acts faithfully on $H_1(X,\Z)$.
\end{lemma}

{\bf Proof.} It suffices to show that the action of any finite subgroup $H\subset \Aut X$ on $H^1(X,\Q)$ is non trivial. Let $Y=X/H$. The quotient map $q:X\to Y$ is a branched covering. Let $n=|H|$, set $k$ to be the number of branch points of $q$ in $Y$, and let $a_1,\ldots, a_k$ be the cardinalities of the preimages of the branch points.  By the Riemann-Hurwitz formula we have
$$\chi(X)=n(\chi(Y)-k)+\sum_{i=1}^k a_i=n\chi(Y)-nk+\sum_{i=1}^k a_i.$$ We have $-nk+\sum_{i=1}^ka_i<0$, so if $\chi(Y)<0$, then $\chi(X)\neq\chi(Y)$. If $\chi(Y)\geq 0$, then we also have $\chi(X)\neq\chi(Y)$ by our assumption on $\chi(X)$. We conclude that $\dim H^1(X,\Q)\neq\dim H^1(Y,\Q)=\dim H^1(X,\Q)^H$, showing that the action of $H$ on $H^1(X,\Q)$ is non trivial.$\clubsuit$

\begin{cor}
The groups $N_\Gamma(G)/G$ and $N_{\PSL_2(\R)}(G)/G$ act faithfully on the abelianisation $G^{ab}$ of $G$.
\end{cor}

{\bf Proof.} Let $E_2, E_3\subset Y(G)$ be the locus of elliptic points of order 2, respectively 3, and set $E=E_2\cup E_3$. If $R$ is a ring and $A$ is a set, we let $R[A]$ be the free $R$-module generated by $A$. We have an $N_{\PSL_2(\R)}(G)/G$-equivariant exact sequence
$$0\to \Z[E]\to H_1(Y(G)\backslash E,\Z)\to H_1(Y(G),\Z)\to 0.$$ Since $\chi(Y(G)\backslash E)<0$, the action of $N_{\PSL_2(\R)}(G)/G$ on the middle term is faithful by lemma~\ref{fin_auto_faithful}. The group $H_1(Y(G),\Z)$ is $G^{ab}$ quotiented by the torsion subgroup $G^{ab}_T$ of $G^{ab}$. Since all groups in the sequence are free Abelian, an element of $N_{\PSL_2(\R)}(G)/G$ that acts trivially on $H_1(Y(G),\Z)$, must act non-trivially on $\Z[E]$, hence on $E$ and hence on $G^{ab}_T=\Z/2[E_2]\oplus \Z/3[E_3]$.$\clubsuit$
%The elliptic points of $Y(G)$ are one to one with the conjugacy classes of finite subgroups: take the preimage in $\Hh$ and then the stabilisers of all elements. The conjugacy classes of finite subgroups are one to one with finite order generators in any given system of free generators obtained from a fundamental domain as in Theorem~\ref{main1}: each such subgroup is conjugate to the group generated by a generator because of the fundamental domain, and for Abelianisation reasons there can not be more congugacy classes. Finally, the torsion of $G^{ab}$ is generated by the images of the torsion elements in any system of free generators obtained from a fundamental domain as in Theorem~\ref{main1}.

\section{Algorithms for congruence subgroups}\label{classical}

The results of the previous sections allow one to construct the bipartite cuboid graph corresponding to a subgroup $G\subset\Gamma$ (and hence, a fundamental domain for the action of $G$ on $\Hh$).
%In general the number of operations this procedure takes is roughly $(\Gamma:G)^2$.
%For general $G$ we do not know a procedure that would work much faster than the naive method sketched in the introduction. However, 
In order to implement this procedure one needs to describe the set $G\backslash \Gamma$ and the right action of $\Gamma$ on it. In this section we do this for the classical congruence subgroups $\Gamma_0(N),\Gamma^0(N),\Gamma_1(N),\Gamma^1(N)$, and $\Gamma(N)$. The idea is as follows.
%if $G$ is a level $N$ congruence subgroup, which we will assume in this section, there is an easier way to construct the bipartite cuboid graph, assuming one can find a quick way to write down a set of right coset representatives of $\SL_2(\Z/N)/\{\pm I\}$ modulo the image of $G$ in $\SL_2(\Z/N)/\{\pm I\}$.% In this section we give the method for the classical congruence subgroups $\Gamma_0(N),\Gamma^0(N),\Gamma_1(N),\Gamma^1(N)$ and $\Gamma(N)$.
%The results of the previous section allow one to construct a fundamental domain and a free system of generators of $G\subset\Gamma$ so that the number of operations is roughly the index $(\Gamma:G)$, assuming one can find a quick way to write down a set of representatives of right conjugacy classes of $\Gamma$ modulo $G$. For the classical congruence subgroups $\Gamma_0(N),\Gamma_1(N)$ and $\Gamma(N)$ the number of operations is $O(N),O(N^2)$ and $O(N^3)$ times something which grows more slowly than any positive power of $N$.

Let $G\subset\Gamma$ be a level $N$ congruence subgroup, i.e.\ $G\supset\Gamma(N)$. For a subgroup $H\subset\Gamma$ set ${}_N H$ to be the image of $H$ in $\SL_2(\Z/N)/\{\pm I\}$. Construct a graph ${}_N\Xc(G)$ exactly in the same way as $\Xc(G)$ in the section~\ref{grfun}, but replacing $\Gamma, G,G_0$ and $G_1$ by ${}_N\Gamma, {}_NG,{}_NG_0$ and ${}_NG_1$ respectively. Then the right $\Gamma$-sets $G\backslash\Gamma$ and ${}_NG\backslash{}_N\Gamma$ are isomorphic. So the graphs $\Xc(G)$ and ${}_N\Xc(G)$ will also be isomorphic, and the isomorphism will respect all the additional structure (partition of the vertices into two types, cyclic order at the trivalent vertices, the distinguished edges). So for level $N$ congruence subgroups it suffices to work with elements of ${}_N\Gamma$ instead of $\Gamma$.

%In this section we illustrate this simplified algorithm in the case when $G$ is one of the classical congruence subgroups $\Gamma_0(N),\Gamma^0(N),\Gamma_1(N),\Gamma^1(N),\Gamma(N)$.
\iffalse
\label{injective}Set $\mathbb{X}=\{\pm I\}\backslash(\Z/N)^2$. (Here we view $(\Z/N)^2$ as the set of columns of length $2$ with entries in $\Z/N$, and $\{\pm I\}\subset\GL_2(\Z/N)$ acts on $(\Z/N)^2$ by multiplication on the left.) There is a natural injective mapping ${}_N\Gamma\to\mathbb{X}^2$ which takes $A\in{}_N\Gamma$ to the couple formed by the images in $\mathbb{X}$ of the first and the second columns of a lift of $A$ to $\SL_2(\Z/N)$. So it is convenient to store the elements of ${}_N\Gamma$ as couples $(x,y)$ with $x,y\in\mathbb{X}$.
\fi

\smallskip

{\bf Conventions.} Recall that a sequence $(a_N)_{N\geq 1}$ of real numbers is said to be $O(P\log N)$ if there is a polynomial $p$ such that for all $N\geq 1$ we have $|a_N|\leq |p(\log N)|$. 

By an \emph{(elementary) operation} we mean integer addition, multiplication, and division with remainder.

Suppose a group $H$ acts on a set $X$, and let $X'$ be a subset of $X$ which intersects each orbit exactly once. Elements of $X'$ will be called {\it orbit representatives}. An algorithm that constructs, given $x\in X$, the element $x'\in X'$ that belongs to the same $H$-orbit as $x$ will be called a {\it reduction procedure}, cf.~remark~\ref{red}. 

Unless stated otherwise, all rings in this section will be assumed associative, commutative, with identity. We denote the group of invertible elements of a ring $R$ by $R^*$.

\begin{remark} 
As we mentioned in the introduction, in~\cite{nie_parent} Z.~Nie and C.~X.~Parent give a different procedure for constructing coset representatives for $G\backslash\Gamma$ when $G=\Gamma_0(N),\Gamma^0(N),\Gamma_1(N),\Gamma^1(N)$, or $\Gamma(N)$.
\end{remark}

\subsection{\texorpdfstring{Lists of representatives for $\Gamma_0(N)$ and $\Gamma^0(N)$}{Lists of representatives for gamma naught}}\label{repg0}
Recall the definition of the projective line $\mathbb{P}^1(\Z/N)$ over $\Z/N$. Set ${}_N\mathbb{Y}$ to be the set of couples $(a,b)\in(\Z/N)^2$ such that $a$ and $b$ are coprime modulo $N$, i.e., for some (hence for any) lifts $\bar a$ and $\bar b$ of $a$ and $b$ to $\Z$ we have $\gcd(\bar a,\bar b,N)=1$. The group $(\Z/N)^*$ acts on ${}_N\mathbb{Y}$ and we set $\mathbb{P}^1(\Z/N)=(\Z/N)^*\backslash {}_N\mathbb{Y}$. We denote the image of $(a,b)\in{}_N\mathbb{Y}$ in $\mathbb{P}^1(\Z/N)$ by $(a:b)$.

If $a,l$ are integers, $l>1$, we denote the image of $a$ in $\Z/l$ by $[a]_l$.

The quotients ${}_N\Gamma_0(N)\backslash{}_N\Gamma$ and ${}_N\Gamma^0(N)\backslash{}_N\Gamma$ are both isomorphic to $\mathbb{P}^1(\Z/N)$. Indeed, ${}_N\Gamma$ acts transitively on the right on $\mathbb{P}^1(\Z/N)$ by the rule
$$(a':b')\cdot\left[\begin{smallmatrix}a&b\\c&d\end{smallmatrix}\right]=(a'a+b'c:a'b+b'd).$$ The stabiliser of $(0:1)$ is ${}_N\Gamma_0(N)$ and the stabiliser of $(1:0)$ is ${}_N\Gamma^0(N)$. So we get isomorphisms $${}_N\Gamma_0(N)\backslash{}_N\Gamma\to\mathbb{P}^1(\Z/N), \mbox{ }
{}_N\Gamma^0(N)\backslash{}_N\Gamma\to\mathbb{P}^1(\Z/N)$$ which take an element $A\in{}_N\Gamma_0(N)\backslash{}_N\Gamma$, respectively $A\in{}_N\Gamma^0(N)\backslash{}_N\Gamma$ to the image in $\mathbb{P}^1(\Z/N)$ of the second, respectively first row of a matrix $A'$ obtained by lifting $A$ first to ${}_N\Gamma$ and then to $\SL_2(\Z/N)$. 
We have the commutative diagrams
\begin{equation}\label{com_diag}
\begin{tikzcd}
\SL_2(\Z/N)\ar[r]\ar[d,"\mbox{{\scriptsize second row}}"'] & {}_N\Gamma_0(N)\backslash{}_N\Gamma\ar[d,"\cong"]\\
{}_N\mathbb{Y}\ar[r] & \mathbb{P}^1(\Z/N),
\end{tikzcd}
\begin{tikzcd}
\SL_2(\Z/N)\ar[r]\ar[d,"\mbox{{\scriptsize first row}}"'] & {}_N\Gamma^0(N)\backslash{}_N\Gamma\ar[d,"\cong"]\\
{}_N\mathbb{Y}\ar[r] & \mathbb{P}^1(\Z/N),
\end{tikzcd}
\end{equation}
in which the horizontal arrows are the quotients maps, and the right vertical arrows have just been described.
%The isomorphism can be obtained as follows: the map from $\SL_2(\Z)$ to $\mathbb{P}^1(\Z/N)$ which takes a matrix $A$ to the image in $\mathbb{P}^1(\Z)$ of the transpose of the first column of $A$ identifies two matrices precisely when one is obtained from the other by multiplication by an element of $\Gamma_0(N)$ on the left.
%from takes an element $\Gamma_0(N)A\in\Gamma_0(N)\backslash\Gamma$ to the image in $\mathbb{P}^1(\Z/N)$ of the first column of a lift of $A$ to $\SL_2(\Z)$.

So to write a list of representatives for $\Gamma_0(N)$ and $\Gamma^0(N)$ it suffices to give a representative $(a,b)\in{}_N \mathbb{Y}$ for any $(a':b')\in\mathbb{P}^1(\Z/N)$.
Note that given a row $([a]_N,[b]_N)\in{}_N\mathbb{Y}$
%$\left(\begin{smallmatrix}[a]_N\\ [b]_N\end{smallmatrix}\right)$
%such that $0,\leq a,b<N,\gcd(a,b,N)=1$,
one can complete it to a matrix $\left(\begin{smallmatrix}[a]_N&[b]_N\\ [c]_N&[d]_N\end{smallmatrix}\right)\in\mathrm{SL}_2(\Z/N)$ in $O(\log N)$ operations (use the Euclidean division algorithm to find integers $c,d$ such that $ad-cd\equiv1\mod N$).

\bigskip

Write $N=p_1^{m_1}\cdot\cdots\cdot p_k^{m_k}$ with $p_1,\ldots, p_k$ distinct primes and all $m_i>0$. By the Chinese remainder theorem there is a ring isomorphism $$\Z/N\to\prod_{i=1}^k \Z/p_i^{m_i},\quad [a]_N\mapsto \left([a]_{p_1^{m_1}},\ldots, [a]_{p_k^{m_k}}\right).$$ Let
$$\psi:\prod_{i=1}^k \Z/p_i^{m_i}\to\Z/N,\quad \left([a_1]_{p_1^{m_1}},\ldots, [a_k]_{p_k^{m_k}}\right)\mapsto \sum_{i=1}^k \left[ a_i\left(\frac{N}{p_i^{m_i}}\right)^{p^{m_i}_i-p^{m_i-1}_i}\right]_N$$ be the inverse isomorphism.

One can extend the definition of the projective line to an arbitrary ring $R$. Namely, we say that $a,b\in R$ are {\it coprime} if together they generate $R$ as an $R$-module, and we define $\mathbb{P}^1(R)$ to be the quotient of the set of coprime couples $(a,b)\in R^2$ by the diagonal action of the group of units of $R$. (The result is the set of $R$-valued points of the scheme $\mathbb{P}^1_R$.) Note that $\SL_2(R)$ acts transitively on the set of coprime couples, hence so does $\PSL_2(R)=\SL_2(R)/\{\pm I\}$ on $\mathbb{P}^1(R)$. Let $\Gamma_0(R)\subset\PSL_2(R)$ be the stabiliser of $(0:1)\in\mathbb{P}^1(R)$.

\begin{lemma}
Let $R_1, R_2$ be rings. We then have a bijection
$$\mathbb{P}^1(R_1\times R_2)\cong \mathbb{P}^1(R_1)\times\mathbb{P}^1(R_2).$$
\end{lemma}

{\bf Proof.} The functors $\GL_2(-), \SL_2(-),\Gamma_0(-)$ etc.\ from rings to groups commute with finite Cartesian products, hence so does the functor $\mathbb{P}^1(-)$ from rings to sets: if $R$ is a ring, we have $\mathbb{P}^1(R)=\Gamma_0(R)\backslash\PSL_2(R)$, functorially in $R$.$\clubsuit$

So there is a bijection $\mathbb{P}^1(\Z/N)\to \prod\mathbb{P}^1(\Z/p_i^{m_i})$ given by \begin{equation}\label{iso2}
\left([a]_N:[b]_N\right)\mapsto\left(([a]_{p_1^{m_1}}:[b]_{p_1^{m_1}}),\ldots,([a]_{p_k^{m_k}}:[b]_{p_k^{m_k}})\right).
\end{equation}
The inverse is
\begin{equation}\label{iso1}
\left(([a_1]_{p_1^{m_1}}:[b_1]_{p_1^{m_1}}),\ldots,([a_k]_{p_k^{m_k}}:[b_k]_{p_k^{m_k}})\right)\mapsto \left(\psi([a_1]_{p_1^{m_1}},\ldots, [a_k]_{p_k^{m_k}} ):\psi([b_1]_{p_1^{m_1}},\ldots, [b_k]_{p_k^{m_k}} )\right).
\end{equation} 
So for our purposes it would suffice to consider the case when $N$ is a power of a prime, $N=p^m$. 

Since $N=p^m$, in each $(\Z/N)^*$-orbit in ${}_N\mathbb{Y}$ there is an element $(a,b)$ with $$a\in\{[0]_N,[1]_N,[p]_N,\ldots,[p^{m-1}]_N\}.$$ For an $a\in\Z/N$ we denote the stabiliser $\subset(\Z/N)^*$ of $a$ by $\Stab(a)$.

\begin{Prop}
We have $\Stab([0]_N)=(\Z/N)^*, \Stab([1]_N)=\{[1]_N\}$. The stabiliser of $[p^i]_N,1\leq i\leq m-1$ is the kernel of the group homomorphism $f^*:(\Z/N)^*\to (\Z/p^{m-i})^*$ induced by the ring homomorphism $f:\Z/N\to \Z/p^{m-i}$.
\end{Prop}

{\bf Proof.} The first two assertions are clear. To prove the third, observe that $[a]_N$ stabilises $[p^i]_N$ iff $[p^i]_N ([a]_N-1)=0$. Now, the elements of $\Z/N$ annihilated by $[p^i]_N$ are precisely the multiples of $[p^{m-i}]_N$, i.e.\ the elements of $\ker f$. So $\Stab([p^i]_N) =1+\ker f$. All elements of $1+\ker f$ are invertible as the elements of $\ker f$ are nilpotent. So we have $1+\ker f=\ker f^*$, and our claim follows.$\clubsuit$

\medskip

\iffalse

If $a\in(\Z/N)^*$ is in the kernel of $f$, then $a=b\cdot [p^{m-i}]_N+1$ for some $b\in\Z/N$, so $a\cdot [p^i]_N=[p^i]_N$.

Conversely, take an $a\in\Stab(p^i)$. Let $\bar a$ be a representative of $a$, $1\leq\bar a< N-1$. Use the Euclidean division to write $\bar a=xp^j+y$ with $x,y\in\Z$ such that $0\leq y\leq p-1$ and either $x=0$ or $x$ is coprime to $p$ and $j>0$.
If $x=0$, then $y=1$ (otherwise $a$ wouldn't stabilise $[p^i]_N$) and $a=[1]_N\in\ker f$.

Suppose $x$ is coprime to $p$ and $1\leq j<m-i$. Since $a$ stabilises $[p^i]_N$, we have $\bar a p^i=xp^{i+j}+yp^i=p^i+lp^m,l\in\Z$, which means that $yp^i-p^i$ belongs to $p^{i+j}\Z$. Moreover, $yp^i-p^i\neq 0$, since $x$ is coprime to $p$. But $|p^i-y\cdot p^i|$ can take values $0,p^i,\ldots,p^i(p-2)$, so $yp^i-p^i$ can not be a non-zero element of $p^{i+j}\Z$.

Suppose $x$ is coprime to $p$ and $j\geq m-i$. Then $\bar a p^i\equiv y p^i\mod N$, which implies $y=1$, and we have $a=b\cdot [p^{m-i}]_N+1\in\ker f$.
$\clubsuit$
\fi

The map $f^*$ is surjective. So $\{[a]_N\mid 1\leq a<p^{m-i},\gcd(a,p)=1\}$ is a set of representatives of $(\Z/N)^*$ modulo the stabiliser of $p^i,1\leq i\leq m-1$.

\begin{Prop}
A set of orbit representatives for the action of $(\Z/N)^*,N=p^m$ on ${}_N\mathbb{Y}$ can be chosen as the set of all $(a,b)$ such that precisely one of the following holds:
\begin{itemize}
\item $a=[0]_N,b=[1]_N$,
\item $a=[1]_N,b\in\Z/N$,
\item $a=[p^i]_N,1\leq i\leq m-1$, $b=\left[\bar b\right]_N,1\leq \bar b< p^{m-i},\gcd(\bar b,p)=1$.
\end{itemize}

Given $(a,b)\in{}_N\mathbb{Y}$, the reduction procedure goes as follows. Lift $a$ and $b$ to integers $\bar a,\bar b,0\leq \bar a,\bar b\leq N-1$. If $\bar a=0$, then the orbit representative is $([0]_N,[1]_N)$. Otherwise write $\bar a=c\cdot\gcd(\bar a,N)$. If $\gcd(\bar a,N)=1$, then the orbit representative of $(a,b)$ is $([1]_N,b[c]_N^{-1})$. If $\gcd(\bar a,N)=p^i,1\leq i\leq m-1$, then the orbit representative of $(a,b)$ is $([p^i]_N,[b']_N)$ where $1\leq b'<p^{m-i}$ is the integer such that the reduction of $b[c]_N^{-1}$ modulo $p^{m-i}$ is $[b']_{p^{m-i}}$. The reduction procedure takes $O(P\log N)$ operations on each $(a,b)\in{}_N \mathbb{Y}$.

Using the isomorphisms~(\ref{iso2}) and~(\ref{iso1}) one obtains a set of orbit representatives and a reduction procedure for the action of $(\Z/N)^*$ on ${}{}_N\mathbb{Y}$ for arbitrary $N$. The reduction procedure takes $O(P\log N)$ operations on each $(a,b)\in{}_N \mathbb{Y}$, assuming the prime factorisation of $N$ is known.

\end{Prop}

$\clubsuit$

\begin{remark}
The indices of the subgroups $\Gamma_0(N),\Gamma_1(N)$ and $\Gamma(N)$ of $\Gamma$ are given by
$$(\Gamma:\Gamma_0(N))=N\prod_{p|N}\left(1+\frac{1}{p}\right), (\Gamma:\Gamma_1(N))=N^2\prod_{p|N}\left(1-\frac{1}{p^2}\right), (\Gamma:\Gamma(N))=N^3\prod_{p|N}\left(1-\frac{1}{p^2}\right)$$
where each product is over the prime divisors of $N$, see e.g.~\cite[1.2]{diashur}. Using
$$\prod_{p|N}\left(1-\frac{1}{p}\right)\geq\frac{1}{e^\gamma\log\log N+\frac{3}{\log\log N}}\mbox{ for }N\geq 3$$ where $\gamma$ is Euler's constant (see e.g.~\cite[Theorem 15]{rosser_schoenfeld}), we conclude that the index of $\Gamma_0(N), \Gamma_1(N)$, respectively $\Gamma(N)$ in $\Gamma$ is bounded below by $N,  O(N^{2-\varepsilon})$ for all $\varepsilon>0$, respectively $O(N^{3-\varepsilon})$ for all $\varepsilon>0$. Any of these expressions exceeds the complexity of factorising $N$ using even the most straightforward algorithms such as trial division. So in the sequel calculating the prime factorisation of $N$ will not affect any of the complexity bounds we give below that involve the index $(\Gamma:G)$ for $G$ in any of the classical series $\Gamma_0(N), \Gamma^0(N), \Gamma_1(N), \Gamma^1(N), \Gamma(N)$.
\end{remark}

Next we consider the left action of ${}_N\Gamma_0(N)$ and ${}_N\Gamma^0(N)$ on ${}_N\Gamma$. We use the diagrams~(\ref{com_diag}).

\begin{Prop}
For any orbit representative $(a,b)$ for the action of $(\Z/N)^*$ on ${}_N\mathbb{Y}$ from the previous proposition construct a matrix $A_{(a,b)}=\left(\begin{smallmatrix}c&d\\a&b\end{smallmatrix}\right)\in\mathrm{SL}_2(\Z/N)$. The set formed by the images of all such matrices in ${}_N\Gamma$ will be the set of orbit representatives for the action of ${}_N\Gamma_0(N)$ on ${}_N\Gamma$.

The reduction procedure is as follows: take an element $A\in{}_N\Gamma$, lift it to a matrix $A'\in\SL_2(\Z/N)$ and apply the reduction procedure from the previous proposition to the second row of $A'$ to obtain a couple $(a,b)$; the image of the matrix $A_{(a,b)}$ in ${}_N\Gamma$ will be the orbit representative of $A$.
The reduction procedure takes $O(P\log N)$ operations.

Note that the element of ${}_N\Gamma_0(N)$ that takes $A_{(a,b)}$ to $A$ is $A\cdot A_{(a,b)}^{-1}$.

The case of $\Gamma^0(N)$ is similar, except given a representative $(a,b)\in{}_N\mathbb{Y}$, we need to construct a matrix $\left(\begin{smallmatrix}a&b\\c&d\end{smallmatrix}\right)\in\mathrm{SL}_2(\Z/N)$, and not $\left(\begin{smallmatrix}c&d\\a&b\end{smallmatrix}\right)$, and in the reduction procedure we replace ``the second row'' by ``the first row''.

The number of operations which are needed to write down all orbit representatives is in both cases $({}_N\Gamma:{}_N\Gamma_0(N))=({}_N\Gamma:{}_N\Gamma^0(N))$ times $O(P\log N)$.
\end{Prop}

$\clubsuit$

\subsection{\texorpdfstring{Lists of representatives for the remaining classical congruence subgroups}{Lists of representatives for other gammas}}\label{repg1g}

We start with a few straightforward observations on how to obtain a list of representatives and a reduction procedure for the action of a subgroup. We will then apply these to cover the remaining classical congruence subgroups. 

Suppose a group $H$ acts freely on a set $X$ on the left and let $\{x_i\}$ be a set of orbit representatives for this action. Suppose $\{h_j\}$ be a set of orbit (=right coset) representatives for the left action of a subgroup $H'\subset H$ on $H$. Then $\{h_j x_i\}$ is a set of orbit representatives for the action of $H'$ on $X$. Suppose we have reduction procedures for the action of $H$ on $X$ and for the left action of $H'$ on $H$ and moreover, suppose that, given $x\in X$ we can not only construct an $x_i$ that belongs to the same $H$-orbit as $x$, but also an element $h\in H'$ such that $hx_i=x$.

Then we can construct a reduction procedure for the action of $H'$ on $X$ as follows. Take an $x\in X$, reduce it modulo $H$ to $x_i$ and find an $h$ in $H$ such that $hx_i=x$. Apply the reduction procedure to $h$ to obtain $h=h'h_j$ with $h'\in H'$, i.e.\ $x=h'h_jx_i$. Note that we obtain not only an orbit representative $h_jx_i$ of $x$, but also an element $h'$ that takes the orbit representative to $x$.

\bigskip

\iffalse
The set of all $\left[\begin{smallmatrix}a&0\\0&a^{-1}\end{smallmatrix}\right]$ such that $a=[\bar a]_N,\gcd(\bar a,N)=1,1\leq \bar a\leq \frac{N}{2}$ is a set of (both left and right) coset representatives of both ${}_N\Gamma_0(N)$ modulo ${}_N\Gamma_1(N)$ and ${}_N\Gamma^0(N)$ modulo ${}_N\Gamma^1(N)$. The reduction procedure both for the action of ${}_N\Gamma_1(N)$ on ${}_N\Gamma_0(N)$ and for the action of ${}_N\Gamma^1(N)$ on ${}_N\Gamma^0(N)$ consists in setting the non-diagonal element equal to zero. The above observation gives us lists of orbit representatives and reduction procedures for the left actions of ${}_N\Gamma_1(N)$ and ${}_N\Gamma^1(N)$ on ${}_N\Gamma$.

As above, the number of operations necessary to write the sets of coset representative is in both cases $({}_N\Gamma:{}_N\Gamma_1(N))=({}_N\Gamma:{}_N\Gamma^1(N))$ times $O(P\log N)$, and the reduction procedures take $O(P\log N)$ operations.

\bigskip

In a similar way, passing from ${}_N\Gamma_1(N)$ to ${}_N\Gamma(N)$ (the identity subgroup) we can list the elements of ${}_N\Gamma$ in $|{}_N\Gamma| O(P\log N)$ operations. We don't need a reduction procedure for the action of ${}_N\Gamma(N)$ on~${}_N\Gamma$.%, since each element is a coset representative modulo the identity subgroup.

\bigskip
\fi

Here is a family of subgroups of $\Gamma$ that includes $\Gamma_0(N),\Gamma_1(N)$, and $\Gamma(N)$. Let $S\subset(\Z/N)^*$ be a subgroup, and let $l$ be a divisor of $N$. Define $\Gamma(N,l;S)\subset\Gamma$ to be the subgroup of all $\left[\begin{smallmatrix}a & b \\c& d\end{smallmatrix}\right]\in\Gamma$ such that $[a]_N,[d]_N\in S, b\equiv 0\mod l, c\equiv 0\mod N$. We have $\Gamma_0(N)=\Gamma(N,1;(\Z/N)^*), \Gamma_1(N)=\Gamma(N,1;e)$, and $\Gamma(N)=\Gamma(N,N;e)$ where $e=\{[1]_N\}$. The group $\Gamma(N,l;S)$ does not change if we replace $S$ with the subgroup of $(\Z/N)^*$ generated by $S$ and $[-1]_N$. Each $\Gamma(N,l;S)\subset \Gamma_0(N)$. 

We want to give a set of representatives and a reduction procedure for the left action of $\Gamma(N,l;S)$ on~$\Gamma_0(N)$. By the above remarks, it suffices to do the same for the left action of ${}_N\Gamma(N,l;S)$ on ${}_N\Gamma_0(N)$. Assume that for the action of $S$ on $(\Z/N)^*$ we have a set of representatives $L$ and a reduction procedure. Without restricting generality we may (and will) also assume that $[-1]_N\in S$.%, and hence that $L$ only contains elements $[\bar a]_N$ with $1\leq \bar a\leq \frac{N}{2}$.

%Note: ${}_N\Gamma(N,1;S)$ is normal in ${}_N\Gamma_0(N)$.
Let us first consider the case $l=1$. To get a set of representatives for the left action of ${}_N\Gamma(N,1;S)$ on ${}_N\Gamma_0(N)$ one can take the set of all $\left[\begin{smallmatrix}a&0\\0&a^{-1}\end{smallmatrix}\right]$ such that $a\in L$. The reduction procedure first sets the non-diagonal element equal zero, then reduces modulo $S$ any lift of the top left element from $(\Z/N)^*/\{\pm [1]_N\}$ to $(\Z/N)^*$. To pass from $\Gamma(N,1;S)$ to $\Gamma(N,l;S)$ for an arbitrary $l|N$ we use the bijection $$\Gamma(N,l;S)\backslash \Gamma(N,1;S)\cong \Gamma(N,l;\pm e)\backslash \Gamma(N,1;\pm e)$$ where $\pm e=\{\pm [1]_N\}$. The quotient set on the right is in fact a group isomorphic to $\Z/l$.
%Note: $\Gamma(N,l;S)$ is in general not normal in $\Gamma(N,1;S)$. For example, if $l=N$ and $S$ is everything, then $\Gamma(N,l;S)$ is the group of elements of $\Gamma$ that become diagonal mod $N$.

For the transpose of $\Gamma(N,l;S)$ the algorithms are completely similar. In particular, this takes care of the group $\Gamma^1(N)$, which is the last remaining classical congruence subgroup that we needed to consider.

\smallskip

For $G$ in any of the classical series $\Gamma_0(N)), \Gamma^0(N), \Gamma_1(N), \Gamma^1(N), \Gamma(N)$, the number of operations needed to write the sets of coset representative is $(\Gamma:G)$ times $O(P\log N)$, and the reduction procedures take $O(P\log N)$ operations.

\subsection{Constructing the graph}\label{constr_the_graph}

In sections~\ref{repg0}-\ref{repg1g} we saw how to construct, for $G$ equal to one of the classical congruence subgroups $\Gamma_0(N), \Gamma^0(N), \Gamma_1(N), \Gamma^1(N), \Gamma(N)$, the set $G\backslash \Gamma$ and the right action of $\Gamma$ on it. We will now explain how to use this data to obtain, for an arbitrary finite index subgroup $G\subset\Gamma$, the graph $\Xc(G)$, additionally assuming that $G\backslash\Gamma$ is well ordered (equivalently, totally ordered). Observe that whenever a group acts on a well-ordered set, one can easily construct the quotient set -- just take the smallest element in each orbit.

To construct $\Xc(G)$ we need to specify the following data:
\begin{enumerate}
\item the sets of vertices of type~(0) and (1);
\item for each vertex $x$, the set $C(x)$ of the vertices $x$ is connected with, and for each $x'\in C(x)$, the number of edges (which is 1 or 2) that join $x$ and $x'$;
\item for each trivalent vertex $x$ of type~(1), the cyclic order of the edges that meet at $x$;
\item the distinguished edge.
\end{enumerate}

The set of vertices of type~(0), respectively of type~(1) is the set of $G_0$-, respectively $G_1$-orbits in $G\backslash\Gamma$, represented each by its smallest element. This gives us item 1 from the list. For item 2, we take a vertex say of type~(1), $G_1(x)=\{x,x', x''\}\subset G\backslash\Gamma$ where $x$ is the smallest element of $G_1(x)$, and $x'=x\cdot \left[\begin{smallmatrix}0&1\\-1&1\end{smallmatrix}\right], x''=x\cdot \left[\begin{smallmatrix}0&1\\-1&1\end{smallmatrix}\right]^2$. The $G_0$-orbits $G_0(x), G_0(x'), G_0(x'')$ are the vertices of type~(0) that $G_1(x)$ is connected with. For $y\in G_1(x)$, the number of edges that connect $G_1(x)$ with $G_0(y)$ is the cardinality of the intersection $G_1(x)\cap G_0(y)\subset G\backslash\Gamma$.

Note that if $x,x',x''$ as above are pairwise distinct, they are cyclically ordered, and the cyclic order depends only on the cyclic ordering of $G_1$, and not on the choice of the well ordering of $G\backslash\Gamma$. This gives us item 3. The distinguished edge is the coset $\mathbf{r}_G=G\in G\backslash \Gamma$. The vertices joined by $\mathbf{r}_G$ are the $G_0$- and $G_1$-orbits of $\mathbf{r}_G$.

Suppose it takes no more than $C>0$ operations to compute the action of a given element of $\Gamma$ on a given element of $G\backslash\Gamma$. For each of the classical congruence subgroup series we can take $C=O(P\log N)$. The procedure we have just described requires $C\cdot O((\Gamma:G))$ operations. Combining this with sections~\ref{repg0}-\ref{repg1g} we get
\begin{Prop}\label{constructing_graph}
For $G$ in any of the classical series $\Gamma_0(N), \Gamma^0(N), \Gamma_1(N), \Gamma^1(N)$, or $\Gamma(N)$, the graph $\Xc(G)$ can be constructed in $$(\Gamma:G)O(P\log N)$$ operations.
\end{Prop}

$\clubsuit$

\subsection{Constructing the polygon}\label{constr_the_poly}

Here we explain how one can construct a special polygon starting from the information contained in the list in section~\ref{constr_the_graph}. In particular, we assume that we have lists of the vertices of types~(0) and~(1) of $\Xc(G)$. We denote these lists $\mathrm{V}_0$ and $\mathrm{V}_1$ respectively and further assume that they are ordered. We will identify a vertex with its number in the corresponding list. We will construct the set of edges of a special polygon $P$ for $G$, and the corresponding involution and independent set of generators.

At each step of the algorithm we will have the following data:

\begin{itemize}
\item Ordered lists $L_0,L_1$ and $\mathcal{G}$ of the same length containing respectively vertices of type~(0), vertices of type~(1) and elements of $\Gamma$. The length of these lists will vary but will never exceed $|\mathrm{V}_0|$. At each step the list $L_0$ contains the vertices of type~(0) that lead to the parts of the graph $\Xc(G)$ that we have not yet investigated. If $V$ is the $j$-th element of $L_0$, then the $j$-th element $V'$ of $L_1$ will be called the \emph{predecessor} of $V$. The $j$-th element of $\mathcal{G}$ is the map $\in\Gamma$ that takes $\ii$ to $V$ and $\rr$ to $V'$.  Some elements of $L_0$ are marked ``to be removed''.
%\item a point $\jmath(x)$ for any element $x$ in $L_0\sqcup L_1$;
\item An auxiliary ordered list $I$ of length $|\mathrm{V}_0|$ which contains the positions of some vertices of type~(0) in $L_0$.% (as we will see, each vertex of type~(0) occurs in $L_0$ at most twice and for some vertices we will have to store the places where they appear for the first time).
\item A marking of each vertex of $\Xc(G)$ of type~(0) as either ``visited'' or ``unvisited''. This information is stored in another ordered list of length $|\mathrm{V}_0|$.%, the initial value being ``unvisited'' for all vertices.
\item A set $E$ of geodesic (half-)intervals in $\Hh$, a fixed point free involution $\sigma$ of $E$ and, for each $s\in E$, a transformation $g_s\in G$ taking $s$ to $\sigma(s)$.
\end{itemize}

When the algorithm terminates, $E$ will be the set of edges of a special polygon for $G$, $\sigma$ will be the corresponding involution, and $g_s$ will be the corresponding independent generators of $G$.

\medskip

For a vertex $V$ of $\Xc(G)$ of type~(1) we set $\mathop{\mathrm{Star}}V$ to be the union of the edges meeting at $V$. Note that $\mathop{\mathrm{Star}}V$ is either a single edge (in which case $V$ is a univalent vertex and corresponds to an elliptic point of $G\setminus\Hh$ of order 3), or a union of three edges. In the latter case the vertices of type~(0) of these edges may be all disjoint or two of these vertices may coincide (in which case $\mathop{\mathrm{Star}}V$ contains a loop). We will call these three graphs respectively a ``segment'', a ``sling shot'' and a ``racket''; they are shown in figures~\ref{fig:1a},~\ref{fig:1e} and~\ref{fig:1c} respectively.

In the sequel, ``adding an item $x$ to an ordered list'' means appending $x$ at the end of the list.%; whenever we set the $j$-th element of $L_0$ equal a vertex $V$, we write $j$ in the entry of $I$ corresponding to $V$.

Here is how we start. All vertices of type~(0) are marked as unvisited, and the lists $L_0,L_1,{\mathcal{G}}$ are empty. Recall that we denote the distinguished vertex of type $i=0,1$ by $V_i^G$. Consider three cases. 

1. $\mathop{\mathrm{Star}}V_1^G$ is a segment. Then we add $V^G_0,V^G_1$ and $\left[\begin{smallmatrix}1&0\\0&1\end{smallmatrix}\right]$ to $L_0, L_1$ and $\mathcal{G}$ respectively, and write 1 in the $V_0^G$-entry of $I$. We add the geodesic half-intervals $s_1=[\rr,0)$ and $s_2=[\rr,\infty)$ to $E$, and we set $\sigma(s_1)=s_2,g_{s_1}=\left[\begin{smallmatrix}0&1\\-1&1\end{smallmatrix}\right]^{-1}.$

2. $\mathop{\mathrm{Star}}V_1^G$ is a sling shot. Let $V'$ and $V''$ be the two vertices of $\mathop{\mathrm{Star}}V_1^G$ of type~(0) other than $V^G_0$; assume the edge connecting $V_1^G$ and $V'$ comes after the distinguished edge in the cyclic ordering of the edges meeting at $V_1^G$. We add $V^G_0,V',V''$ to $L_0$ (in this order) and write 1,2 and 3 in corresponding the entries of $I$. The predecessor of each of these is set to be $V_1^G$, and the corresponding elements of $\mathcal{G}$ are defined to be $\left[\begin{smallmatrix}1&0\\0&1\end{smallmatrix}\right],\left[\begin{smallmatrix}0&1\\-1&1\end{smallmatrix}\right]$ and $\left[\begin{smallmatrix}0&1\\-1&1\end{smallmatrix}\right]^2$ respectively. We mark $V^G_0,V',V''$ as visited.

3. $\mathop{\mathrm{Star}}V_1^G$ is a racket. Let $V$ be the vertex of type~(0) of $\mathop{\mathrm{Star}}V_1^G$ which is at the tip of the handle of the racket. Let $j$ be 0,1 or 2, depending on whether the distinguished edge coincides with the handle or comes before or after it in the cyclic ordering of the edges meeting at $V_1^G$. We add $V, V_1^G$ and $\left[\begin{smallmatrix}0&1\\-1&1\end{smallmatrix}\right]^j$ to $L_0,L_1$ and $\mathcal{G}$ respectively; write 1 in the $V$-entry of $I$. Let $s_0=(0,\infty),s_1=(0,1)$ and $s_2=(1,\infty)$. We add $s_{j+1}$ and $s_{j+2}$ to $E$, and we define $\sigma(s_{j+1})=s_{j+2},g_{s_{j+1}}=\left[\begin{smallmatrix}0&1\\-1&1\end{smallmatrix}\right]^{j-1}\left[\begin{smallmatrix}1&-1\\0&1\end{smallmatrix}\right]\left[\begin{smallmatrix}0&1\\-1&1\end{smallmatrix}\right]^{1-j}$.% We mark the vertex of $\mathop{\mathrm{Star}}V_1^G$ which at the head of the racket as ``visited''.

Now we describe the induction step. Let $j$ be the current length of $L_0$. If $j=0$, we stop, otherwise we continue.

Take the last entry $V$ of $L_0$; let $V'$ and $g$ be the corresponding entries of $L_1$ and $\mathcal{G}$. The first thing to do is to remove $V,V'$ and $g$ from $L_0,L_1$ and $\mathcal{G}$ respectively (of course, after having saved them somewhere for the induction step). If $V$ was marked ``to be removed'', the induction step is complete. Otherwise we proceed as follows.

%1. Suppose $V$ is marked as ``visited''. Then $V$ already occurs at some place, say $l$, in $L_0$; we can find $l$ using the list $I$. Let $g'$ be the $l$-th entry of $\mathcal{G}$. Let $s_0$ be the geodesic connecting 0 and $\infty$ and set $s_1=g\cdot s_0,s_2=g'\cdot s_0$. Add $s_1$ and $s_2$ to the list $E$, set $\sigma(s_1)=s_2$ and $g_{s_1}=g'g^{-1}$. We mark the $l$-th entry of $L_0$ ``to be removed''.
%
1. Suppose $V$ is
% marked as ``unvisited''. If $V$ is
univalent. Then set $s'_1=[\ii,0)$ and $s'_2=[\ii,\infty)$. We add $s_1=g\cdot s'_1$ and $s_2=g\cdot s'_2$ to $E$, define $\sigma(s_1)$ to be $s_2$ and set $g_{s_1}=g\left[\begin{smallmatrix}0&1\\-1&0\end{smallmatrix}\right] g^{-1}$. Note that $\left[\begin{smallmatrix}0&1\\-1&0\end{smallmatrix}\right]$ takes $s_1'$ to $s_2'$ while keeping $\ii$ fixed.

2. Suppose $V$
% is marked as ``unvisited'' and
has valency 2. Let $V''$ be the vertex of type~(1) other than $V'$ that $V$ is connected with by an edge. As above, there are three cases to consider.

2a. $\mathop{\mathrm{Star}}V''$ is a segment. Then let $s_1'=[\mathrm{e}^{\frac{2\pi\ii}{3}},0)$ and $s_2'=[\mathrm{e}^{\frac{2\pi\ii}{3}},\infty)$. We add $s_1=g\cdot s'_1$ and $s_2=g\cdot s'_2$ to $E$, define $\sigma(s_1)$ to be $s_2$ and set $g_{s_1}=g\left[\begin{smallmatrix}1&1\\-1&0\end{smallmatrix}\right] g^{-1}$. Note that $\left[\begin{smallmatrix}1&1\\-1&0\end{smallmatrix}\right]$ takes $s_1'$ to $s_2'$ while keeping $\mathrm{e}^{\frac{2\pi\ii}{3}}$ fixed.

2b. $\mathop{\mathrm{Star}}V''$ is a sling shot.
% Then we mark $V$ as ``visited''.
Let $V_1$ and $V_2$ be the vertices of $\mathop{\mathrm{Star}}V''$ of type~(0) such that the edge connecting $V''$ and $V_1$ comes after the edge connecting $V''$ and $V$ in the cyclic ordering of the edges meeting at $V''$. Let $W\in\{V_1,V_2\}$, and set $\tilde g$ to be $g\left[\begin{smallmatrix}1&-1\\0&1\end{smallmatrix}\right]$ if $W=V_1$ and $g\left[\begin{smallmatrix}1&0\\-1&1\end{smallmatrix}\right]$ if $W=V_2$. There are the following possibilities.

2b$'$. If $W$ is marked as unvisited, then we add it to $L_0$ (recall that $V$ has been deleted). We write the current length of $L_0$ in the $W$-entry of the list $I$. The predecessor of $W$ is set to be $V''$. The last element of $\mathcal{G}$ is set to be $\tilde g$. We mark $W$ as visited.

2b$''$. If $W$ is marked as visited, then it already occurs at some place, say $l$, in $L_0$; we can find $l$ using the list $I$. Let $g'$ be the $l$-th entry of $\mathcal{G}$. Let $s'=(0,\infty)$ and set $s_1=\tilde g\cdot s',s_2=g'\cdot s'$. Add $s_1$ and $s_2$ to $E$, and set $\sigma(s_1)=s_2$ and $g_{s_1}=g'\left[\begin{smallmatrix}0&1\\-1&0\end{smallmatrix}\right]\tilde g^{-1}$. We mark the $l$-th entry of $L_0$ ``to be removed''.

2c. $\mathop{\mathrm{Star}}V''$ is a racket. Let $s'=(0,\infty)$. Set $g_1=g\left[\begin{smallmatrix}1&-1\\0&1\end{smallmatrix}\right]$ and $g_2=g\left[\begin{smallmatrix}1&0\\-1&1\end{smallmatrix}\right]$. We add $s_1=g_1\cdot s'$ and $s_2=g_2\cdot s'$ to $E$ and define $\sigma(s_1)=s_2,g_{s_1}=g_2\left[\begin{smallmatrix}0&1\\-1&0\end{smallmatrix}\right]g_1^{-1}$.

\medskip

The procedure we have just seen takes $O(1)(\Gamma:G)$ operations. Using proposition~\ref{constructing_graph} we get
\begin{Prop}\label{constructing_polygon}
For $G$ in any of the classical series $\Gamma_0(N)), \Gamma^0(N), \Gamma_1(N), \Gamma^1(N), \Gamma(N)$, a fundamental special polygon for the action of $G$ on $\Hh$ can be constructed in $$(\Gamma:G)O(P\log N)$$ operations.
\end{Prop}

$\clubsuit$

\bigskip

\noindent
Alexey Gorinov: Faculty of Mathematics, Higher School of Economics, 6 Usacheva ulitsa, Moscow, Russia 119048, email: \url{agorinov@hse.ru}, \url{gorinov@mccme.ru}.

\noindent
Isaac Kalinkin: email: \url{kalinkin.isaac@icloud.com}.

\end{document}